\documentclass[final]{siamltex}

\RequirePackage{hypernat}%
\RequirePackage{amsmath}%
\RequirePackage{amsxtra}%
\RequirePackage{amstext}%
\RequirePackage{amssymb}%
\RequirePackage{latexsym}%
\RequirePackage{dsfont}%
\RequirePackage{ifthen}%
\RequirePackage[left,mathlines]{lineno}%
\RequirePackage[colorlinks,citecolor=blue,urlcolor=blue,pagebackref]{hyperref}
\RequirePackage{color}%

\definecolor{backcol1}{rgb}{0.95,1,1}%
\definecolor{backcol2}{rgb}{1,0.95,1}%
\definecolor{backcol3}{rgb}{1,1,0.95}%
\definecolor{backcol4}{rgb}{1,0.95,0.95}%
\definecolor{backcol5}{rgb}{0.95,0.95,1}%
\definecolor{backcol6}{rgb}{0.95,1,0.95}%
\definecolor{backcol7}{rgb}{0.95,0.9,0.75}%

\definecolor{hellgelb}{rgb}{1,1,0.85}%
\definecolor{colKeys}{rgb}{0,0,1}%
\definecolor{colIdentifier}{rgb}{0,0,0}%
\definecolor{colComments}{rgb}{1,0,0}%
\definecolor{colString}{rgb}{0,0.5,0}%

\definecolor{red}       {rgb}{0.0,0.0,1.0}        % linkcolor
\definecolor{magenta}   {rgb}{0.0,0.0,1.0}        % url local file color
\definecolor{cyan}      {rgb}{0.0,0.0,1.0}        % url color
\definecolor{green}     {rgb}{0.0,0.4,0.3}        % cite color

\newcommand{\eqdef}{%
\ensuremath{\mathrel{\stackrel{\mathrm{def}}{=}}}}

\begin{document}

\title{Parallelizing Sequential Sweeping on Structured Grids -- Fully Parallel
SOR/ILU preconditioners for Structured n-Diagonal Matrices}

\author{R. Tavakoli \thanks{Department of Material Science and Engineering,
Sharif University of Technology, Tehran, Iran, P.O. Box 11365-9466,
email:
\href{mailto:tav@mehr.sharif.edu}{tav@mehr.sharif.edu}, %
\href{mailto:rohtav@gmail.com}{rohtav@gmail.com}.%
}%
}%

\maketitle

\begin{center}
{{\tiny\today}}\vspace*{5mm}%
 %{{\scriptsize 5 November 2009}}\vspace*{5mm}%
\end{center}

\date{\today}

%===============================================================

\begin{abstract}%
There are variety of computational algorithms need sequential
sweeping; sweeping based on specific order; on a structured grid,
e.g., preconditioning (smoothing) by SOR or ILU methods and solution
of eikonal equation by fast sweeping algorithm. Due to sequential
nature, parallel implementation of these algorithms usually leads to
miss of efficiency; e.g. a significant convergence rate decay.
Therefore, there is an interest to parallelize sequential sweeping
procedures, keeping the efficiency of the original method
simultaneously. This paper goals to parallelize sequential sweeping
algorithms on structured grids, with emphasis on SOR and ILU
preconditioners. The presented method can be accounted as an
overlapping domain decomposition method combined to a multi-frontal
sweeping procedure. The implementation of method in one and two
dimensions are discussed in details. The extension to higher
dimensions and general structured n-diagonal matrices is outlined.
Introducing notion of alternatively block upper-lower triangular
matrices, the convergence theory is established in general cases.
Numerical results on model problems show that, unlike related
alternative parallel methods, the convergence rate and efficiency of
the presented method is close to the original sequential method.
Numerical results also support successful use of the presented
method as a cache efficient solver in sequential computations as
well.%
\end{abstract}%

%===============================================================



%===============================================================
%===============================================================

\section{Introduction}
\label{sec:int}

Structured grids are commonly used in scientific computing for the
purpose of numerical solution of partial differential equations
(PDEs). Despite the unstructured grid based methods which manage
domain complexity in a more consistent manner, in certain
applications using structured grids is still the best choice, mainly
due to preference in terms of consumed memory, efficiency and the
ease of implementation. Moreover, with the invention of immersed
boundary methods with make possible to manage complex geometries on
structured grids, there is recently a significant interest to
application of structured grids in scientific computing. Progress in
power of supercomputer changes the bias in favor of structured grids
too, as they are usually more appropriate for parallel computing.

Some of numerical algorithms on structured grids use sequential
sweeping; grid by grid procedure on a mandatory order. In fact,
using the conventional version of sequential algorithms, it is not
possible to perform computations simultaneously. This limit us to
use full power of multi-processor shared and/or distributed memory
machines.

Since the advent of parallel computers, several parallel versions of
the SOR method have been developed. Most variants of these
algorithms have been developed by using multi-color ordering schemes
or domain decomposition techniques
\cite{adams1982mcs,oa1985lar,adams1986scb,neumann1987cpm,white1987mpi,melhem1987tei,melhem1988mrs,ashcraft1988vif,block1990bcs,ortega1991ocg,ii1993oma,stotland1997opc,xie1999nps}.
In contrast to domain decomposition method, a multi-color SOR method
has usually a faster convergence rate, but is usually more difficult
to solve elliptic boundary value problems in complex domains
\cite{xie2006nbp}. To improve the convergence rate of a domain
decomposition algorithms, a novel mesh partition strategy was
proposed in \cite{xie1999nps} which led to an efficient parallel
SOR. This method has the same asymptotic rate of convergence as the
Red-Black SOR \cite{xie1999nps}. Although, multi-color SOR
algorithms are also usually implemented based on a domain
decomposition strategy, the Xie and Adams's algorithm
\cite{xie1999nps} takes less inter-processor data communication time
than a multi-color SOR. The blocked version of the Xie and Adams's
parallel SOR algorithm \cite{xie1999nps} is suggested in Xie
\cite{xie2006nbp}.  According to presented results, the rate of
convergence and parallel performance of this version is better than
that of the non-blocked version.

The main drawback of the above mentioned parallel versions of SORs
is lower rate of convergence in contrast to original one. This is
mainly due to enforced decoupling between equations to provide some
degree of concurrency in computations. For example in the m-color
SOR method, an increase in the number of colors decreases the rate
of convergence. So a key to achieve an efficient parallel SOR method
is to minimize the degree of decoupling due to parallelization.

Other techniques such as the pipelining of computation and
communication and an optimal schedule of a feasible number of
processors are also applied to develop parallel SOR methods when the
related coefficient matrix is banded
\cite{missirlis1987spi,robert1988csp,eisenstat1989csp,niethammer1989smp,comput1995pgs,amodio1996pis}.
These techniques can isolate parts of the SOR method which can be
implemented in parallel without changing the sequential SOR. The
convergence rate of these methods are same as the original one but
the application is limited to banded-structure matrices.

The present goals to suggest a new method to parallelize sequential
sweeping on structured grids. The specific focus is on the SOR
method which can be regarded as an iterative preconditioner,
smoother or even a linear solver. Without loss of generality, the
later case is taken into account in this study to simplify
arguments. Since there is a direct correspondence between SOR and
ILU preconditioners on structured grids, the presented method can be
used to develop parallel ILU methods as well.

The rest of this paper is organized as follows. Sections
\ref{sec:1d} and \ref{sec:2d} present the implementation of method
in one and two spatial dimensions respectively. The extension to
three spatial dimensions is commented in Section \ref{sec:3d}. he
Considering the connection between SOR and ILU preconditioners on
structured grids, Section \ref{sec:ilu} argues the application of
presented approach to develop ILU(p) preconditioners on structured
grids. The convergence theory is established in Section
\ref{sec:conv_anal}. Section \ref{sec:gextension} extend utility of
the presented method to general structured n-diagonal matrices.
Section \ref{sec:res} is devoted to numerical results on the
convergence and performance of the method. Finally Section
\ref{sec:closing_remark} close this paper by summarizing results and
outlining possible application and extension of the presented
method.

%===============================================================
%===============================================================

\section{Parallel SOR sweeping in one-dimension}%
\label{sec:1d}

Consider the following scaler 1D elliptic PDE:
\begin{equation}%
\label{eq:1d}%
-\frac{d}{dx}\big(\alpha\ \frac{du}{dx}\big) + \beta u = f \
\mathtt{in} \ \in \Omega, \ \mathtt{and} \ u = u_0(x) \ \mathtt{on} \ \Gamma%
\end{equation}%
where $\Omega = [0, L_x]$, $\Gamma = \partial \Omega$, $\alpha:
\Omega \rightarrow \mathbb{R}^+$ is $C^1$, $\beta \in \mathbb{R}^+
\cup \{0\}$. Also, it is assumed field variables $u$, $u_0$ and $f$
posses the sufficient regularity required by the applied numerical
method to ensure the existence and uniqueness of the discrete
solution. Assume the spatial domain is discrtized into
an $n+1$ intervals with grid points %
$x_i= x_{i-1} + \delta x_{i-1},\ i=1, \ldots, n$  and $x_0 = 0, \
x_{n+1} = L_x$, where $\delta x_i = x_{i+1}-x_{i}$.

Without loss of generality; for the ease of presentation; we
discrtize \eqref{eq:1d} with a second order finite difference method
as follows:%
\begin{equation}%
\label{eq:d1d}%
-c_i u_{i-1} + b_i u_i - a_i u_{i+1} = f_i, %
\qquad i= 1, \ldots, n.
\end{equation}
where subscript $i$ denotes the value of field variable at $x=x_i$
(e.g. $u_i = u(x_i)$), %
\[a_i = \frac{2\alpha_{i+1/2}}{\delta x_i (\delta x_i+\delta x_{i-1})}, \quad  %
c_i = \frac{2\alpha_{i-1/2}}{\delta x_{i-1} (\delta x_i+\delta
x_{i-1})}, \quad
b_i = a_i + c_i + \beta, \]%
\[
\alpha_{i+1/2} =
\frac{2\alpha_i\alpha_{i+1}}{\alpha_i+\alpha_{i+1}}, \quad
\alpha_{i-1/2}
=\frac{2\alpha_i\alpha_{i-1}}{\alpha_i+\alpha_{i-1}}\].
Using boundary conditions $u_0 = u_0(x_0)$ and $u_{n+1} =
u_0(x_{n+1})$, it is possible to write above equations in the
following matrix form, %
\begin{equation}%
\label{eq:1d:mat}
 A \mathbf{u} = \mathbf{f},
\end{equation}
\[A= \left[ \begin{array}{ccccc}%
b_1&-a_1&&&\\
-c_2&b_2&-a_2& & \\
&-c_3&\ddots&\ddots&\\
&&\ddots & \ddots & -a_{n-1}\\
&&&-c_n&b_n %
\end{array} \right]_{n\times n},\]
\\ %
where
$\mathbf{u} = [u_1, u_2, \ldots, u_n]^T$ and %
$\mathbf{f} = [c_1u_0 + f_1, f_2, \ldots, f_{n-1}, f_n + a_n
u_{n+1}]^T$. It is evident that matrix $A$ is irreducibly diagonally
dominant which ensure the convergence of stationary iterative
methods, like SOR. The SOR method using the natural row-wise
ordering (left to right sweeping) generates the following sequence
of iterations from a
given initial guess $u_i^{(0)}$, %
\begin{equation}%
\label{eq:1d:sorl} %
u_i^{(k+1)}=(1-\omega_{L})u_i^{(k)}+ %
\frac{\omega_{L}}{b_i}%
(c_i u_{i-1}^{(k+1)} + a_i u_{i+1}^{(k)} + f_i), %
\qquad i=1, \ldots, n,
\end{equation}%
where $\omega_{L} \in (0,2)$ is left to right (LR) sweeping
relaxation factor and superscript $(k)$ denotes the iteration
number. The other alternative SOR method is resulted by reversing
the sweeping direction (right to left sweeping):%
\begin{equation}%
\label{eq:1d:sorr}%
u_i^{(k+1)}=(1-\omega_{R})u_i^{(k)}+ %
\frac{\omega_{R}}{b_i}(c_i u_{i-1}^{(k)}+a_i u_{i+1}^{(k+1)} + f_i),
\qquad i=n, \ldots, 1.
\end{equation}%
where $\omega_{R} \in (0,2)$ is the right to left (RL) sweeping
relaxation factor. The application of \eqref{eq:1d:sorl} and
\eqref{eq:1d:sorr} alternatively during each two consecutive
iterations leads to the classical symmetric SOR (SSOR) method.  It
is clear that these procedures are sequential, in the sense that
update of grid point $i$ should be performed after its left (or
right) neighbor. In fact, it is not possible to perform relaxation
of more than one grid point at each time in the original version of
SOR method.

The goal of parallel implementation of SOR method is to provide
possibility of concurrent computation, simultaneously keep the
convergence rate of original SOR method as much as possible. The
concept of domain decomposition is used here to parallelize the SOR
method. Assume spatial domain $\Omega$ is decomposed into $p$ number
of non-overlapping sub-domains, i.e., $\Omega = \cup_{i=1}^p
\Omega_i$, and $\Omega_i \cap_{i\neq j} \Omega_j = \emptyset, \
\mathtt{for} \ i=1, \dots, p$ (splitting of grid points).

\begin{figure}[ht]%
\centering{%
\includegraphics[width=12.cm]{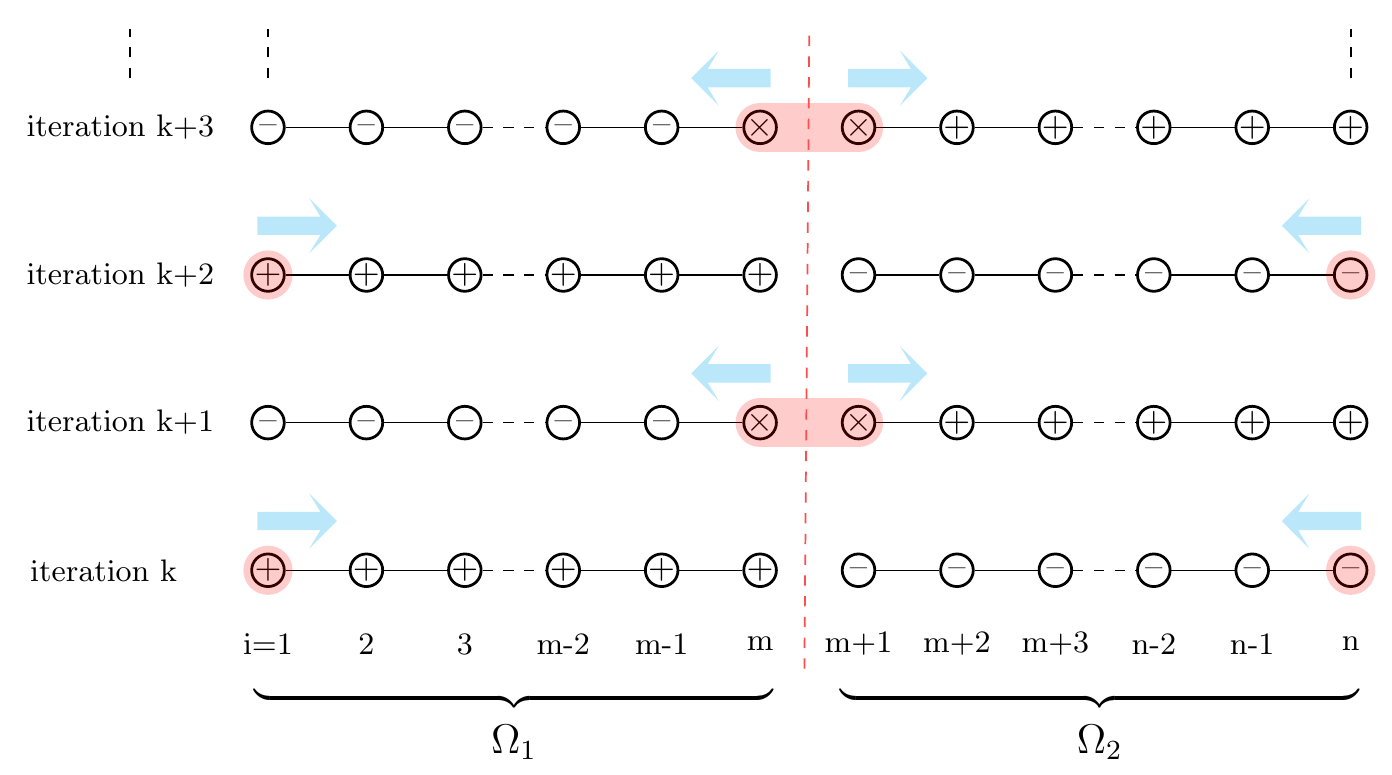}%
}%
\caption{The diagram of the parallel implementation of the presented
SOR method for two sub-domains.}%
\label{fig:fig1}%
\end{figure}

For the purpose of convenience, we first describe our method for two
(approximately equal) sub-domains. As illustrated in figure
\ref{fig:fig1}, discrtized spatial domain is decomposed into
sub-domains $\Omega_1$ and $\Omega_2$, where $\Omega_1 = [0,x_m]$
and $\Omega_2 = [x_{m+1},L]$. If \eqref{eq:1d:sorl} is used in
$\Omega_1$ and the computation begins from the left boundary, and
simultaneously, \eqref{eq:1d:sorr} is used in $\Omega_2$ and the
calculation is performed from the right boundary toward left
direction, then the computation of sub-domains will be independent
during the first iteration. In our algorithm during the next
iteration \eqref{eq:1d:sorr} and \eqref{eq:1d:sorl} are used
respectively in $\Omega_1$ and $\Omega_2$; and the computations are
started from grid points $x_m$ and $x_{m+1}$ in the corresponding
sub-domains. Therefore, we have the following equation at grid point
$x_m$,
\begin{equation}%
\label{eq:1d:sorcl}%
u_m^{(k+1)}-\frac{\omega_{R}a_m}{b_m}u_{m+1}^{(k+1)}=
(1-\omega_{R})u_m^{(k)}+
\frac{\omega_{R}}{b_m}(c_mu_{m-1}^{(k)}+f_m),
\end{equation}
while the following equation should be used at grid point $x_{m+1}$,
\begin{equation}%
\label{eq:1d:sorcr}%
u_{m+1}^{(k+1)}
-\frac{\omega_{L}c_{m+1}}{b_{m+1}}u_{m}^{(k+1)}=(1-\omega_{L})u_{m+1}^{(k)}+
\frac{\omega_{L}}{b_{m+1}}(a_{m+1} u_{m+2}^{(k)}+f_{m+1}).
\end{equation}
\ref{eq:1d:sorcl} and \eqref{eq:1d:sorcr} can be written in the
following matrix form:%
\begin{equation}%
\label{eq:1d:sorc} %
\left[ \begin{array}{cc}
1&-\frac{a_m\omega_{R}}{b_m}\\
-\frac{c_{m+1}\omega_{L}}{b_{m+1}}&1\\
\end{array} \right]
{\left[ \begin{array}{l}
u_{m}^{(k+1)}\\
u_{m+1}^{(k+1)}
\end{array} \right]}
= {\left[ \begin{array}{l}
r_m\\
r_{m+1}
\end{array} \right]},
\end{equation}
\[r_m = (1-\omega_{R})u_{m}^{(k)}+\frac{\omega_{R}}{b_{m}}(c_mu_{m-1}^{(k)} + f_m),\]
\[ r_{m+1}= (1-\omega_{L})u_{m+1}^{(k)}+ \frac{\omega_{L}}{b_{m+1}}(a_{m+1}u_{m+2}^{(k)}+ f_{m+1}). \]
The above $2\times2$ system of equation can be easily inverted to
compute $u_{m+1}^{(k)}$ and $u_{m+1}^{(k+1)}$. Computing the new
values of grid points $x_m$ and $x_{m+1}$ based on the mentioned
procedure, the computations for the remaining grid points will be
performed trivially according to the corresponding update formulae.
These sequence of computation is repeated for each pair of
consecutive iterations until the convergence criteria is satisfied.
The flow chart of this procedure is displayed in figure
\ref{fig:fig1}. In this plot $\oplus$ and $\ominus$ symbols are used
to denote the application of \eqref{eq:1d:sorl} and
\eqref{eq:1d:sorr} respectively; also symbol $\otimes$ is used to
denote the application of \eqref{eq:1d:sorc}.

\begin{figure}[ht]%
\centering{%
\includegraphics[width=13.cm]{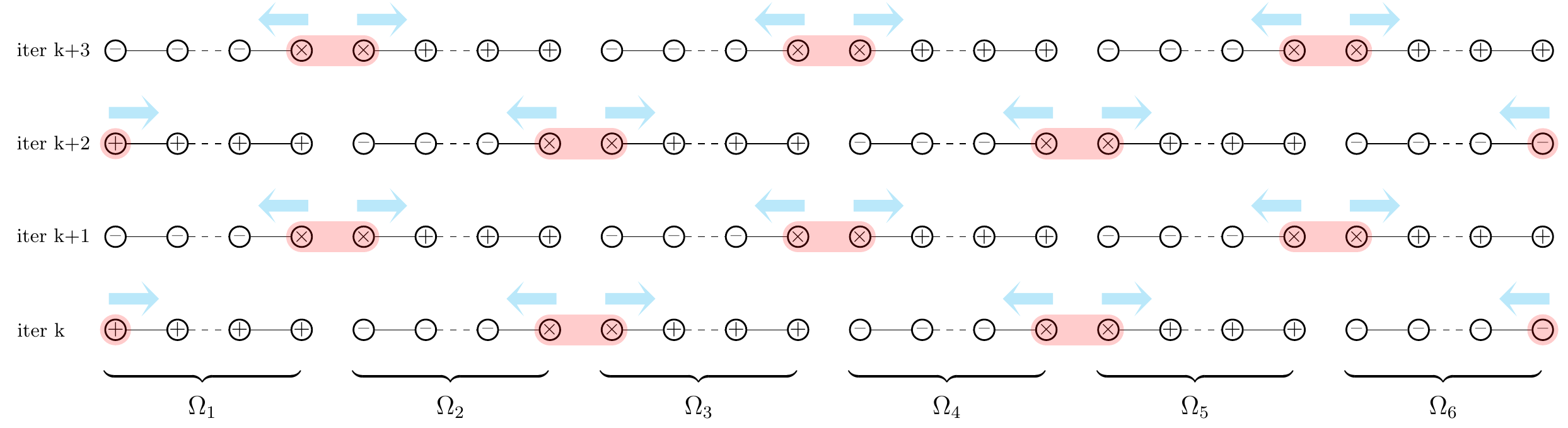}%
}%
\caption{The diagram of the parallel implementation of the presented
SOR method in one-dimension for multiple number of sub-domains.}%
\label{fig:fig2}%
\end{figure}

The extension of this algorithm to more than two sub-domains is a
straightforward job which can be described as follow (cf. figure
\ref{fig:fig2}):
\begin{enumerate}
\item[(i)]
Division of the spatial computational domain into desired number of
sub-domains, based on desired load balancing constraints.%
\item[(ii)]
Determination of the sweeping direction for each sub-domain.
Sweeping direction of each sub-domain should be in opposite
direction of its neighbors. For example we use LR direction for odd
sub-domains and RL direction for even sub-domains. These sweeping
directions are inverted after each iteration.%
\item[(iii)]
Updating the starting node of each sub-domain with
\eqref{eq:1d:sorc} and the remained nodes with either of
\eqref{eq:1d:sorl} or \eqref{eq:1d:sorr}, based on the corresponding
sweeping direction. If the starting node be located at the physical
boundaries, the prescribed boundary value
will be used there (there is no need to use \eqref{eq:1d:sorc} in this case).%
\end{enumerate}

Now let's to look at the properties of the presented parallel
algorithm and its connection to other domain decomposition methods.
It is obvious that this algorithm is an overlapping domain
decomposition method which changes the overlap regions area
alternatively. As the method has some similarities to overlapping
Schwarz algorithm, lets to mention some differences of our algorithm
with the Schwarz method. In contrast to the Schwarz method the size
of overlap regions here are changes alternatively. Sub-system of
equations are only
 solved exactly at overlapped regions (of course alternatively);
 and within internal regions of sub-domains, equations are solved
 approximately (by one-pass SOR iteration). In Schwarz algorithm
 sub-domain boundaries are treated as virtual boundary conditions
 of Dirichlet, Neumann, Robin or mixture of these types. When
 relaxation factors are equal to unity (Gauss–Seidel iterations),
 our numerical treatment for sub-domain boundaries are equivalent
  to Neumann-type (flux) boundary conditions, in the other cases
  our internal boundary conditions are equivalent to a mixture
 of Dirichlet and Neumann ones. In the later case, the relaxation
 factors $\omega_L$ and $\omega_R$ plays role of blending factors.
 In fact our method tries to keep the partial coupling between sub-domain
 during iterations, in expense of a negligible computational cost.
 It is worth mentioning that due to alteration of sweeping
 directions we expect good smoothing properties of the presented
 method. Later we will observe that how such enforced local coupling leads
 to competitive convergence rate of the presented method in contrast
 to that of original SOR.

Note that the procedure for higher order stencils or other numerical
methods are conceptually the same. We avoid further remarks in this
regard to save the space.

%===============================================================
%===============================================================

\section{Parallel SOR sweeping in two-dimensions}%
\label{sec:2d}

In this section we shall extend the presented parallel sweeping
algorithm for two-dimensional problems. We consider the following
scaler 2D elliptic PDE:
\begin{equation}%
\label{eq:2d}%
-\nabla \cdot (\alpha\ \nabla u) + \beta u = f \
\mathtt{in} \ \in \Omega, \ \mathtt{and} \ u = u_0(x) \ \mathtt{on} \ \Gamma%
\end{equation}%
where $\Omega \in \mathbb{R}^2$, $\Gamma = \partial \Omega$,
$\alpha: \Omega \rightarrow \mathbb{R}^{2\times2}$ is $C^1$ diagonal
matrix with diagonal entries $\alpha_x, \alpha_y \in \mathbb{R}^+$,
$\beta \in \mathbb{R}^+ \cup \{0\}$. For the sake of convenience, it
is assumed that $\Omega = [0, L_x] \times [0, L_y]$. Same as
previous, sufficient regularities are assumed for the field
variables.

Suppose that $\Omega$ is divided into $(n+2)$ and $(m+2)$ Cartesian
grids along $x$ and $y$ directions respectively. The grid points $(
x_i , y_j )$ are given by $x_i = x_{i-1} + \delta x_{i-1},\ i = 1,
\ldots, n$ and $y_j = y_{j-1} + \delta y_{j-1},\ j = 1, \ldots, m$;
where $x_0 = y_0 = 0$, $x_{n+1} = L_x$, $y_{m+1} = L_y$, $\delta x_i
= x_{i+1}-x_{i}$ and $\delta y_j = y_{j+1}-y_{j}$. Same as previous,
the finite difference approximation of \eqref{eq:2d} is used here to
describe the presented method (with second order central
differencing scheme). Higher order finite difference methods as well
as other numerical methods (like FVM, FEM) can also be employed. The
only requirement of our method is discretization of PDE on a
structured grid. We use $w_{i, j}$ to denote the finite difference
approximations of $w ( x_i , y_j )$. The approximation of
\eqref{eq:2d} with the mentioned scheme has the following form,
\begin{equation}%
\label{eq:d2d}%
- a_{i,j}^S u_{i,j-1} - a_{i,j}^W u_{i-1,j} + a_{i,j}^P u_{i,j} -
a_{i,j}^E u_{i+1,j}- a_{i,j}^N u_{i,j+1}= f_{i,j},
\end{equation}%
for  $i = 1, \ldots, n$ and $j = 1, \ldots, m$; where $a_{i,j}^P =
a^S_{i,j} + a^W_{i,j} + a^E_{i,j} + a^N_{i,j} + \beta$,
\[
a_{i,j}^E = \frac{2\alpha_{x  i+1/2, j}}{\delta x_{i  } (\delta
x_i+\delta x_{i-1})}, \quad %
a_{i,j}^W = \frac{2\alpha_{x  i-1/2, j}}{\delta x_{i-1} (\delta
x_i+\delta
x_{i-1})}, \]%
\[
a_{i,j}^N = \frac{2\alpha_{y  i, j+1/2}}{\delta y_{j  } (\delta
y_j+\delta y_{j-1})}, \quad %
a_{i,j}^S = \frac{2\alpha_{y  i, j-1/2}}{\delta y_{j-1} (\delta
y_j+\delta y_{j-1})},
\]
and the mid grid point coefficient are computed by harmonic
averaging, e.g.,
\[ \alpha_{x  i+1/2, j} = \frac{2\alpha_{x  i, j} \alpha_{x i+1, j}}
{\alpha_{x  i, j} +\alpha_{x i+1, j}}, \quad
\alpha_{y  i, j+1/2} = \frac{2\alpha_{y  i, j} \alpha_{y i, j+1}}
{\alpha_{y  i, j} +\alpha_{y i, j+1} }.\]%
Using the corresponding boundary conditions, %
\[u_{0, j} = u_0(0,y_j), \quad %
  u_{n+1, j} = u_0(x_{n+1},y_j),\quad %
  j=1, \ldots, m.\]%
\[u_{i, 0} = u_0(x_i,0), \quad %
u_{i, m+1} = u_0(x_i, y_{m+1}), \quad %
i = 1, \ldots, n.\]%
Same as previous, it is possible to write above equations in a
matrix form. However, the presented method is matrix-free and there
is no need to write equations explicitly in a matrix form at this
point.

\begin{figure}[ht]%
\centering{%
\includegraphics[width=10.cm]{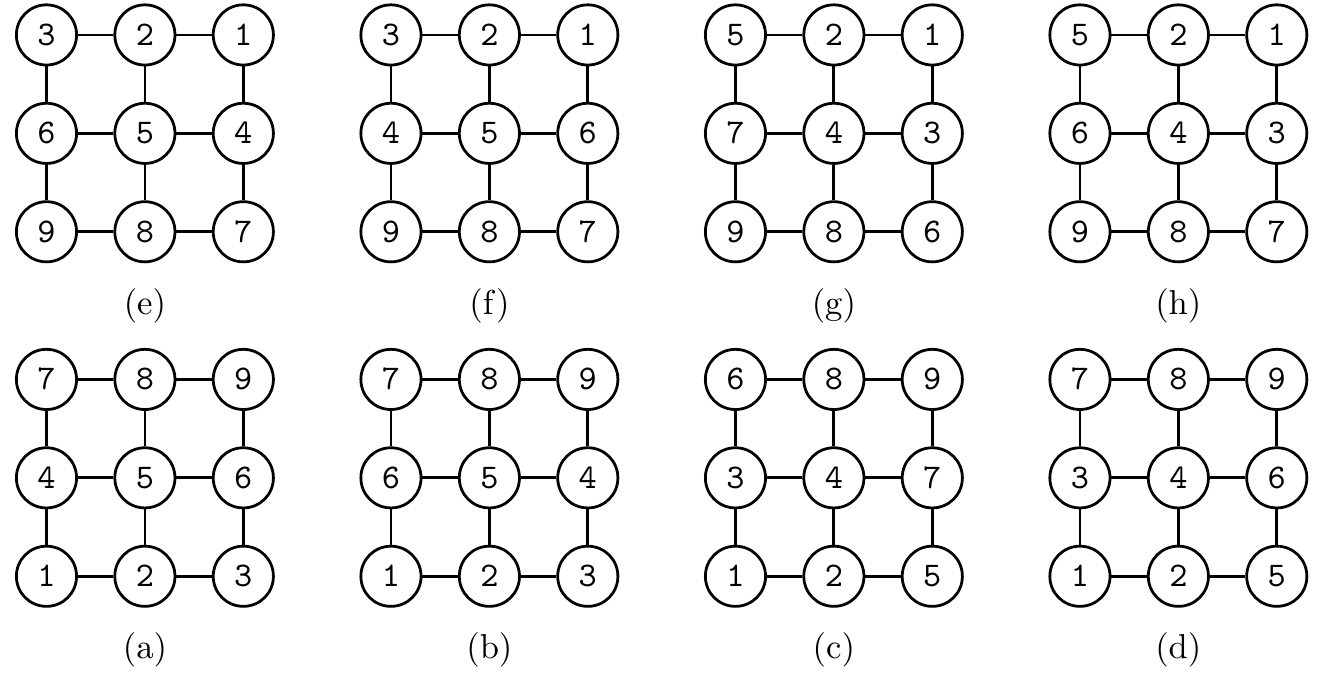}%
}%
\caption{Possible sweeping strategies for  SW-to-NE (bottom row),
NE-to-SW (top row) directions on a $3\times 3$ structured grid:
order of updating procedure in natural row-wise (a, e), symmetric
row-wise (b, f) and frontal (c, d, g, h) sweeping strategy.}%
\label{fig:swp2d}%
\end{figure}

There are variety of sweeping direction to solve \eqref{eq:d2d} with
SOR method. Considering a rectangular computational domain, it is
possible to start sweeping procedure from one corner toward its
corresponding opposite corner. Therefore, there is four sweep
directions. Using geographic directions (east, west, north, south),
it is possible to show these directions by SW-to-NE, NE-to-SW,
SE-to-NW and NW-to-SE. There are various sweeping strategy for each
of mentioned directions. These strategies are shown in figure
\ref{fig:swp2d} for SW-to-NE, NE-to-SW directions on a $3\times3$
structured grid. Every sweeping strategy should preserve causality;
should be performed along a causal direction. The causality here
means that order of update should be such that new values of at
least one half of neighboring nodes of the current grid should be
known prior to its update. It is possible to use distinct relaxation
factors along each sweeping direction.  We show relaxation parameter
of each direction by subscript according to its starting corner,
i.e., $w_{SW}$, $w_{NE}$, $w_{SE}$, $w_{NW}$.

The SOR updating (from iteration $k$ to $k+1$) formulae for grid
point $(i,j)$ has one of the following formulae,

\begin{align}%
\label{eq:sor2d_sw}%
u_{i,j}^{(k+1)} &= (1-\omega_{SW}) u_{i,j}^{(k)} + %
\frac{\omega_{SW}}{a^P_{i,j}}( %
a^S_{i,j} u_{i,j-1}^{(k+1)} + %
a^W_{i,j} u_{i-1,j}^{(k+1)} + %
a^E_{i,j} u_{i+1,j}^{(k)}   + %
a^N_{i,j} u_{i,j+1}^{(k)}   + f_{i,j}) \\%
\label{eq:sor2d_ne}%
u_{i,j}^{(k+1)} &= (1-\omega_{NE}) u_{i,j}^{(k)} + %
\frac{\omega_{NE}}{a^P_{i,j}}( %
a^S_{i,j} u_{i,j-1}^{(k)} + %
a^W_{i,j} u_{i-1,j}^{(k)} + %
a^E_{i,j} u_{i+1,j}^{(k+1)}   + %
a^N_{i,j} u_{i,j+1}^{(k+1)}   + f_{i,j}) \\%
\label{eq:sor2d_se}%
u_{i,j}^{(k+1)} &= (1-\omega_{SE}) u_{i,j}^{(k)} + %
\frac{\omega_{SE}}{a^P_{i,j}}( %
a^S_{i,j} u_{i,j-1}^{(k+1)} + %
a^W_{i,j} u_{i-1,j}^{(k)} + %
a^E_{i,j} u_{i+1,j}^{(k+1)}   + %
a^N_{i,j} u_{i,j+1}^{(k)}   + f_{i,j}) \\%
\label{eq:sor2d_nw}%
u_{i,j}^{(k+1)} &= (1-\omega_{NW}) u_{i,j}^{(k)} + %
\frac{\omega_{NW}}{a^P_{i,j}}( %
a^S_{i,j} u_{i,j-1}^{(k)} + %
a^W_{i,j} u_{i-1,j}^{(k+1)} + %
a^E_{i,j} u_{i+1,j}^{(k)}   + %
a^N_{i,j} u_{i,j+1}^{(k+1)}   + f_{i,j})%
\end{align}

\begin{figure}[ht]%
\centering{%
\includegraphics[width=4.cm]{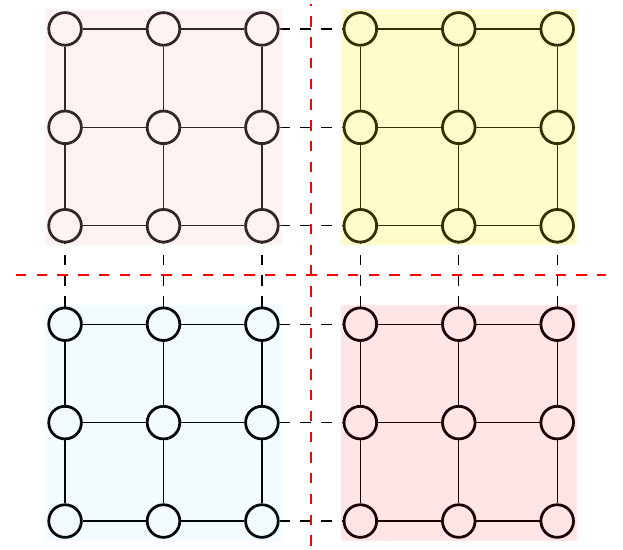}%
}%
\caption{Decomposition of a $6\times6$ Cartesian mesh into
$2\times2$  array of sub-domains.}%
\label{fig:2d_decomp}%
\end{figure}

It is clear that due to required causality condition, the above
procedures are sequential in nature. To parallelize this sequential
procedure, we first split the computational domain into $p$ (almost)
equally-size non-overlapping sub-domains (based on load balancing
constraints). To describe our algorithm, we consider a Cartesian
decomposition in this section, i.e., $p = p_x \times p_y$ where
$p_x$ along $p_y$ are number of sub-domains along $x$ and $y$
spatial axes respectively. For the sake of convenience, it is
assumed that $\Omega$ is decomposed into a $2\times 2$ array of
sub-domains (cf. figure \ref{fig:2d_decomp}).

As it is implied in previous, a key to achieve a parallel SOR with a
little convergence rate decay is to minimize degree of decoupling
due to parallelizing, as much as possible. In the original SOR when
a point is updated the remaining points sense its new value due to
causality of sweeping strategy.

\begin{figure}[ht]%
\centering{%
\includegraphics[width=13.cm]{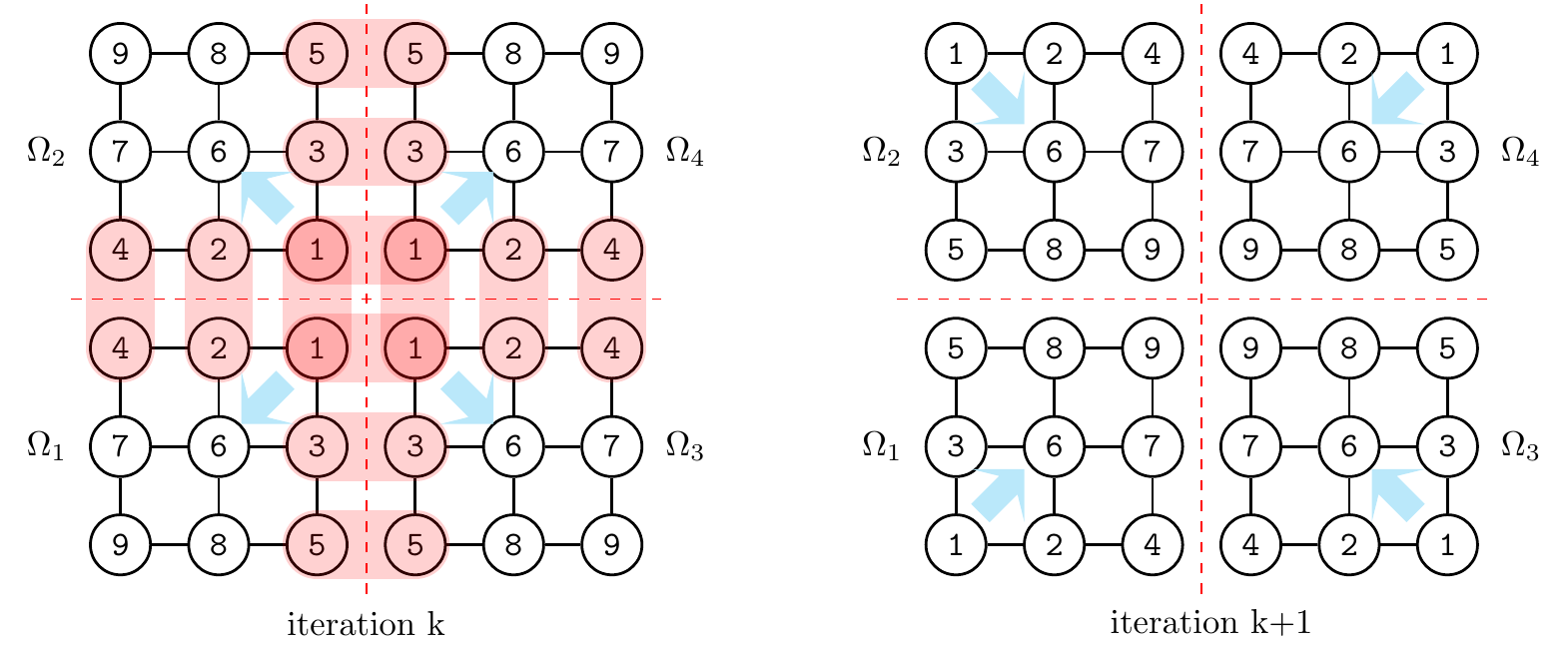}%
}%
\caption{Decomposition of a $6\times6$ Cartesian mesh into
$2\times2$ array of sub-domains: sweeping directions and order of
update are shown during two consecutive iterations, coupling between
nodes are shown by shading.}%
\label{fig:split2d}%
\end{figure}%

To keep the causality of sweeping partially, a multi-frontal
strategy is employed here for the purpose of parallelization.
Consider the mentioned $2\times2$ Cartesian splitting of the
computational domain. According to figure \ref{fig:split2d}, we use
frontal sweeping along directions, $NE$, $NW$, $SE$ and $SW$ for
sub-domains 1 though 4 respectively. The order of updating is shown
in plots. It is clear that to update node \#1 of each sub-domain, we
need two unknown values (new values) from two neighbor sub-domains.
A shaded region in figure show this coupling in iteration $k$-th.
Therefore a local $4\times 4$ system of equation should be solved
for this purpose. Assuming node \#1 of $\Omega_1$  (at iteration
$k$-th, left side of figure) is denoted by $(i,j)$ and using global
index, this local system of equation has the following form
(application of \eqref{eq:sor2d_ne}, \eqref{eq:sor2d_nw},
\eqref{eq:sor2d_se} and \eqref{eq:sor2d_sw} for $\Omega_1$ though
$\Omega_4$ respectively),
\begin{align}%
\label{eq:2d:sor4b4} %
\nonumber%
& \left[ %
\begin{array}{cccc}%
1 & %
-\frac{\omega_{NE} a^E_{i,j}}{a^P_{i,j}} & %
-\frac{\omega_{NE} a^N_{i,j}}{a^P_{i,j}} & %
0 %
\\\\%
-\frac{\omega_{NW} a^W_{i+1,j}}{a^P_{i+1,j}} & %
1 & %
0 & %
-\frac{\omega_{NW} a^N_{i+1,j}}{a^P_{i+1,j}}  %
\\\\%
-\frac{\omega_{SE} a^S_{i,j+1}}{a^P_{i,j+1}} & %
0 & %
1 & %
-\frac{\omega_{SE} a^E_{i,j+1}}{a^P_{i,j+1}} %
\\\\ %
0 & %
-\frac{\omega_{SW} a^S_{i+1,j+1}}{a^P_{i+1,j+1}} & %
-\frac{\omega_{SW} a^W_{i+1,j+1}}{a^P_{i+1,j+1}} & %
1  %
\end{array}%
\right]%
{\left[ %
\begin{array}{l}
u_{i,j}^{(k+1)}\\\\ %
u_{i+1,j}^{(k+1)}\\\\ %
u_{i,j+1}^{(k+1)}\\\\ %
u_{i+1,j+1}^{(k+1)} %
\end{array}%
\right]} =
\\\nonumber\\ %
& \left[ %
\begin{array}{l}%
(1-\omega_{NE}) u_{i,j}^{(k)} + \frac{\omega_{NE} ( a_{i,j}^S
u_{i,j-1}^{(k)} + a_{i,j}^W u_{i-1,j}^{(k)} + f_{i,j})
}{a_{i,j}^P}\\
(1-\omega_{NW}) u_{i+1,j}^{(k)} + \frac{\omega_{NW} (a_{i+1,j}^S
u_{i+1,j-1}^{(k)} + a_{i+1,j}^E u_{i+2,j}^{(k)} + f_{i+1,j})
}{a_{i,j+1}^P}
 \\
(1-\omega_{SE}) u_{i,j+1}^{(k)} + \frac{\omega_{SE} (a_{i,j+1}^W
u_{i-1,j+1}^{(k)} + a_{i,j+1}^N u_{i,j+2}^{(k)} + f_{i,j+1})
}{a_{i,j+1}^P}  \\ %
(1-\omega_{SW}) u_{i+1,j+1}^{(k)} + \frac{\omega_{SW} (a_{i+1,j+1}^E
u_{i+2,j+1}^{(k)} + a_{i+1,j+1}^N u_{i+1,j+2}^{(k)} + f_{i+1,j+1})
}{a_{i+1,j+1}^P} %
\end{array}%
\right]%
\end{align}

This system of equation is non-singular and can be easily inverted.
After update of the corner node (node 1), we should update nodes on
(virtual) boundaries of sub-domains. As it is shown in figure
\ref{fig:split2d}), node \#2 of each sub-domain is coupled to node
\#2 of a  neighbor sub-domain. Assuming node \#1 of $\Omega_1$ is
denoted by $(i,j)$ and using global index, using \eqref{eq:sor2d_ne}
and \eqref{eq:sor2d_se}the following $2\times2$ system of equations
should be solved to compute new values of node \#2 in $\Omega_1$ and
$\Omega_3$, in fact $(i-1,j)$ and $(i-1, j+1)$ (cf.
\ref{fig:split2d}),
\begin{align}%
\label{eq:2d:sor2b2} %
\nonumber%
&\left[ %
\begin{array}{cc}%
1 & %
-\frac{\omega_{NE} a^N_{i-1,j}}{a^P_{i-1,j}} \\\\ %
-\frac{\omega_{SE} a^S_{i-1,j+1}}{a^P_{i-1,j+1}} & %
1 \\\\%
\end{array}%
\right]%
{\left[ %
\begin{array}{l}
u_{i-1,j}^{(k+1)}\\\\ %
u_{i-1,j+1}^{(k+1)} %
\end{array}%
\right]} =
\\\nonumber\\
&{\left[ %
\begin{array}{l}
(1-\omega_{NE}) u_{i-1,j}^{(k)} + \frac{\omega_{NE} (a^S_{i-1,j}
u_{i-1,j-1}^{(k)} + a^W_{i-1,j} u_{i-2,j}^{(k)} + a^E_{i-1,j}
u_{i,j}^{(k+1)})
}{a^P_{i-1,j}}%
\\\\ %
(1-\omega_{SE}) u_{i-1,j}^{(k)} + \frac{\omega_{SE} (a^W_{i-1,j+1}
u_{i-2,j+1}^{(k)} + a^E_{i-1,j+1}u_{i,j+1}^{(k+1)} + a^N_{i-1,j+1}
u_{i-1,j+2}^{(k)})
}{a^P_{i-1,j+1}}%
\end{array}%
\right]}
\end{align}%
Note that the new values of $u_{i,j}^{(k+1)}$ and
$u_{i,j+1}^{(k+1)}$ are located in right hand side of
\eqref{eq:2d:sor2b2} as they are known after solution of
\eqref{eq:2d:sor4b4}. In the same manner other remaining boundary
nodes will be updated. Then, the internal nodes will be updated
using a classic SOR method along the corresponding sweeping
direction.

\begin{figure}[ht]%
\centering{\includegraphics[width=13.cm]{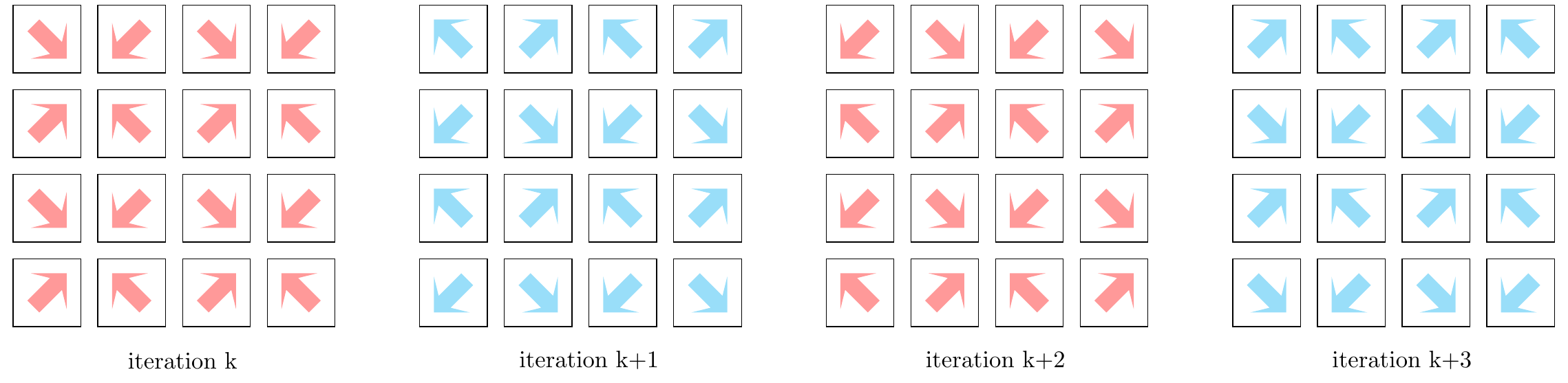}}%
\caption{Alternation of sweeping directions during 4 consecutive
iterations for $4\times4$ array of sub-domains.}%
\label{fig:2d4b4}%
\end{figure}

The next sweep (iteration $(k+1)$-th) for the mentioned
decomposition will be performed along the corresponding opposite
directions, i.e., $SW$, $SE$, $NW$ and $NE$ for sub-domains 1 though
4 respectively (cf. figure \ref{fig:split2d}). In this specific
example there is no need to solution some local system of equations
for this sweep; but in general one may usually needs to solve such
systems at coupling points of sub-domains which are determined based
on sweeping directions. Figure \ref{fig:2d4b4} shows the alteration
of sweeping directions during each four consecutive iterations for a
$4\times4$ array of sub-domains. In the next following sub-sections
some implementation issues of our algorithm will be discussed in
details.

As a final remark in this section it should be mentioned that the
extension of our algorithm for nine-point or higher order central
stencils is straightforward. As the extend of computational molecule
is increased in a higher order method, the size of local systems to
be solved will be also increased accordingly. Following the same
strategy, it would be possible to apply the presented method for
other numerical methods, e.g., FVM and FEM methods, on a structured
grid. Note that every structured grid has a one-to-one
correspondence with with a Cartesian grid. Therefore, the presented
method can be directly used for non-Cartesian structured grids too.
To save the space, further details in this regards are ignored in
this study.

%=========================================================================
%=========================================================================

\subsubsection{Domain decomposition}%
\label{sec:2d:dd}

The first aspect of our parallel algorithm is to divide the
computational grid into parts. This is how the full computational
task is divided among the various processors. Each processor works
on a portion of the grid and anytime a processor needs information
from another one a message will be passed (in message passing
protocol; no message passing is needed regarding to shared-memory
machines).

To achieve a good parallel performance, we like to have an optimal
load balancing and the least communication between processors.
Consider the load balancing first. Assuming a homogeneous parallel
machine, one would like that each processor does the same amount of
work during every iteration. That way, each processor is always
working and not sitting idle. It makes sense to partition the grid
such that each processor gets an equal (or nearly equal) number of
nodes to work on. The second criterion, needs to minimize the number
of edge cuts and the degree of adjacency for each sub-domain (number
of neighbors processor). So the domain decomposition could be
converted to a constrained optimization problem in which the
communication cost to be minimized under the load balancing
constraint. In the case of regular Cartesian grid, this optimization
problem can be easily solved within a negligible cost.

The cost of communication is composed from two types of elapsed
time, one proportional to the communicated data size (function of
network communication speed) and one due to the network latency
which is independent from the data size. Therefore, dealing with a
high latency networks, we may prefer to decrease the degree of
adjacency in our communication graph; in expense of a higher size of
communicated data size. On the other hand, dealing with low latency
networks, we may neglect the effect of network latency. It is worth
mentioning that, using an overlapped communication and computation
strategy, it will be possible to decrease cost of communication.

The above discussion is correct in general for every parallel
algorithms. However, one should consider other criteria when the
convergence rate of the parallel algorithm is a function domain
decomposition topology (number of sub-domains and their geometries).

In the present study increasing number of sub-domains along each
spatial direction is equivalent to decreasing the degree of coupling
in that direction. Rough speaking, it usually decreases the
convergence rate. However, this conclusion is not correct in
general. In some cases, a limited degree of decoupling not only
decreases the convergence rate but also improves it. As an example
consider the solution of Poisson equation in a square domain. Also,
assume that the center of domain is the center of symmetry for the
exact solution. Now, let's to decompose the domain into $2\times2$
equal sub-domains. Due to mentioned symmetry, the parallel algorithm
performs superior than the original SOR method (as it is more
consistent to mentioned symmetry in actual solution). This example
implies that decoupling does not always decrease the convergence
rate. Therefore, one may expect that an algorithm with partial
coupling (like that of ours) performs even superior than the
original SOR in practice.

%=========================================================================
%=========================================================================

\subsubsection{Sweeping directions}%
\label{sec:2d:sw}

After the domain decomposition, it is essential to determine the
sweeping directions for each sub-domain during 4 consecutive
iterations. For this purpose it is sufficient to use the following
rules (lie 1D version of algorithm):\\%
\begin{itemize}%
\item[-]  Along each spatial direction the sweeping direction of a sub-domain
should be in opposite direction of its neighbor sub-domains (cf.
figure \ref{fig:2d4b4}).
\item[-] The sweeping directions (along all spatial directions)
should be reversed during every pair of iterations (cf. figure
\ref{fig:2d4b4}).
\item[-] After each pair of iteration, the sweeping will be
performed along a new diagonal, e.g., after two consecutive
iterations along diagonal $SW$-$NE$, the direction of sweeping for
the next two iteration will be changed to $SE$-$NW$ diagonal (cf.
figure \ref{fig:2d4b4}).\\

\end{itemize}%
Using the above mentioned rules, the sweeping directions are
generated automatically without any complexity, e.g., figure
\ref{fig:2d4b4} is produced applying the above mentioned rules
within a nested loop using a few lines of
Tikz\footnote{\href{http://www.texample.net/tikz/}{www.texample.net/tikz/}}
scripts.

%=========================================================================
%=========================================================================

\subsubsection{Data structure}%
\label{sec:2d:ds}

In the sequential SOR a two-dimensional array is sufficient to store
the global data. This needs a nested loop to sweep the hole domain.
Since in the presented parallel SOR method, we do not have a regular
order of updating for sub-domains boundary nodes, it is preferable
to store order of update for these nodes. To achieve better
performance every data (mentioned indexes and field variables) are
stored in a one-dimensional array (sometimes called as
vectorization). In fact in this way we have a one-to-one mapping
from 2D indexing to 1D indexing, e.g., node $(i, j)$ is stored at
location $i+ j*(n+2)$ in the one-dimensional array. In this way, the
location of nodes $(i+1, j)$ and $(i, j+1)$ will be the location of
node $(i, j)$ plus 1 and $(n+2)$ values respectively.

\begin{figure}[ht]%
\centering{%
\includegraphics[width=10.cm]{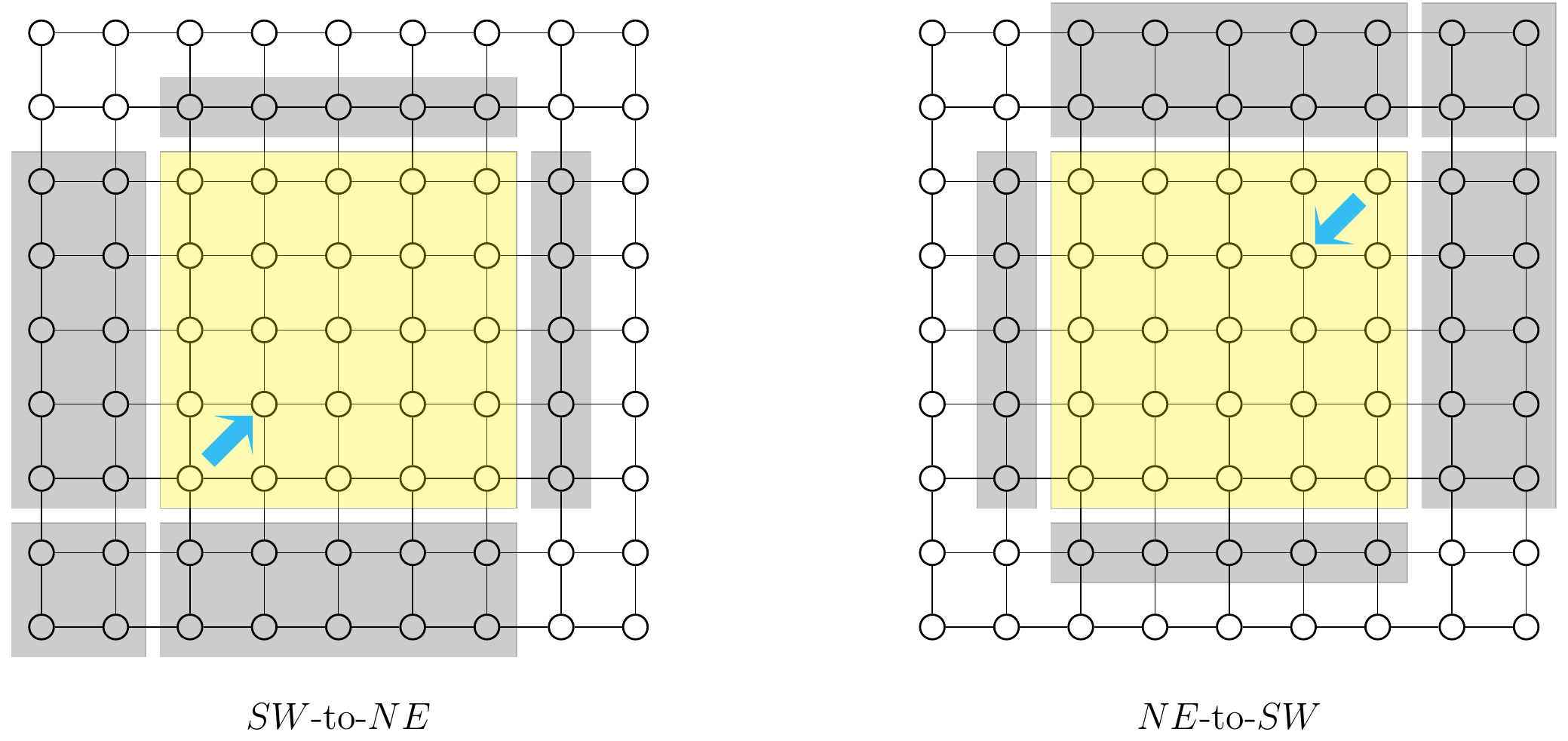}%
}%
\caption{Schematic of data communication during frontal sweeping
$SW$ to $NE$ (left) and $NE$ to $SW$ (right) for a sub-domain with 8
active neighbors: yellow shaded regions contain native nodes and
gray shaded regions include communicated nodes.}%
\label{fig:2dcomm}%
\end{figure}

To update nodes which are located at boundaries between neighbor
sub-domains, the data from the corresponding neighbor sub-domains is
required. Therefore, we add two rows/columns of halo points around
each sub-domain which are renewed (via communication) prior to a new
iteration (in the case of higher order schemes, the width of halo
region should be increased accordingly). In the present
implementation every sub-domain has at most 8 neighbors, i.e., 4
first-degree neighbors (contacted via an edge), 4 second-degree
neighbors (contacted via a corner). Figure \ref{fig:2dcomm} shows
schematically the communicated data during $SW$ to $NE$ and $NE$ to
$SW$ sweeps for a sub-domain with 8 active neighbors. In this
figure, yellow-shaded regions contain native nodes and gray-shaded
regions includes required data which should be communicated. It is
clear that it is possible to overlap a portion of communications
with computations to improve the performance.

%

%===============================================================
%===============================================================

%=====================================================================
%=====================================================================
%=====================================================================
%=====================================================================
%=====================================================================

\section{Extension to three-dimensions}%
\label{sec:3d}

Following the geometric procedure discussed in the previous section,
it is straightforward to extend the presented method to
three-dimensional (3D) cases. The main ingredients of such
extensions are briefly addressed here.

In 3D case, there are eight alternatives frontal sweeping
directions. Adding front and back directions to mentioned geographic
directions (east, west, north, south), these eight sweep directions
are as follows: $BSW$-to-$FNE$, $FNE$-to-$BSW$, $BNE$-to-$FSW$,
$FSW$-to-$BNE$, $BNW$-to-$FSE$, $FSE$-to-$BNW$, $BSE$-to-$FNW$ and
$FNW$-to-$BSE$. There are eight relaxation parameters accordingly,
$\omega_{BSW}$, $\omega_{FNE}$, $\omega_{BNE}$, $\omega_{FSW}$,
$\omega_{BNW}$, $\omega_{FSE}$, $\omega_{BSE}$ and $\omega_{FNW}$.
After the domain decomposition, the sweeping directions are
determined for each domain for each eight consecutive iterations
using mentioned rules in \ref{sec:2d:sw}. Using a Cartesian
decomposition topology, each sub-domain will has at most 26
neighbors, i.e., 6 first-degree neighbors (contacted via a face), 12
second-degree neighbors (contacted via  an edge), 8 third-degree
neighbors (contacted via a corner).

Same as 2D version of algorithm, some local system of equations are
needed to be solved at sub-domains coupling nodes. Now, consider a
sub-domain with 26 active neighbor. At starting corner, a $8\times8$
system of linear equations should be solved. Note that this system
is not dense (includes 32 none-zero entries). Then, along each of
the three coupling edges, a sort of $4\times4$ local system of
linear equations should be solved in an appropriate order. In the
same manner, on each of the three corresponding coupling faces a
sort of $2\times2$ local systems of equations should be solved based
on an appropriate update order. Updating coupling nodes, the
remaining nodes are updated using classic SOR method along the
corresponding sweep direction.

%=====================================================================
%=====================================================================
%=====================================================================
%=====================================================================
%=====================================================================

\section{Parallel incomplete LU preconditioners}%
\label{sec:ilu}%

The SOR method does not usually used as an iterative solver in
practice. But it is commonly used as a smoother in multi-grid
methods or as a preconditioner in Krylov subspace methods. In these
cases, a few iterations of SOR method is applied during each step of
preconditioning (or smoothing). Another kind of preconditioners
which are extensively used in Krylov subspace methods are incomplete
LU (ILU) preconditioners. However due to sequential nature of
incomplete factorization and also forward/backward subsituation
procedure, application of ILU methods in parallel is challenge.

In this section, we show that the presented parallel method can be
equivalently used to develop parallel ILU preconditioners. A simple
observation shows that single pass of a symmetric Gauss-Seidel
method on a n-diagonal matrix (precisely a matrix includes only a
few non-zero diagonals with symmetric sparsity pattern) is
equivalent to ILU(0) preconditioning (cf. section 10.3 of
\cite{saad2003ims}). Therefore, the presented parallel algorithm can
be used as a parallel ILU(0) preconditioner as well. Similarly,
ILU(p) preconditioning for such structured matrices is equivalent to
application of our method for a higher order method (cf. section
10.3 of \cite{saad2003ims}) which make sense to use our method to
parallelize ILU(p) preconditioners as well.

Later, we will show that our method can be applied for n-diagonal
matrices raised from discretization of PDEs on structured grids.
Therefore, the presented strategy can be also regarded as a parallel
ILU(p) preconditioner for such n-diagonal matrices.

%===============================================================
%===============================================================

\section{Convergence analysis}%
\label{sec:conv_anal}

In this section we first proof the convergence of the presented
parallel algorithm in one- and two- spatial dimensions. Following
the same line of proof, the extension of theory to higher dimension
will be straightforward.

\subsection{Mono-dimensional convergence analysis}%
\label{sec:conv_anal_1d}

To simplify the analysis, it is worth to write the presented method
in an algebraic form. Without loss of generality, assume that the
spatial domain is decomposed into $p$ number of sub-domains where
$p$ is an even integer and all sub-domains include $l$ number of
grid points (note that boundary nodes are not considered here, i.e.,
$n=pl$). Also, assume that the sweeping direction of the first
sub-domain is $LR$ at the starting iteration. We can decompose
matrix $A$
 in \eqref{eq:1d:mat} as follows, %
\begin{equation}%
\label{eq:ca1d:decomp}
 A=  G_L + \beta I + G_R,
\end{equation}
where $I$ denotes an $n\times n$ identity matrix here,
\begin{equation}%
\nonumber
G_L = \left[ %
\begin{array}{lllll}%
L_1&&&&  \\%
&R_2&&&  \\%
&&\ddots&& \\%
&&&L_{p-1}&  \\%
&&&&R_p %
\end{array}%
\right],\quad%
G_R = \left[ %
\begin{array}{lllll}%
R_1&&&&  \\%
&L_2&&&  \\%
&&\ddots&& \\%
&&&R_{p-1}&  \\%
&&&&L_p  %
\end{array}%
\right],%
\end{equation}%
for $q=2, 3, \ldots, p$,
\[
L_q = \left[ %
\begin{array}{llllll}%
-c_{(q-1)l+1}&c_{(q-1)l+1}&&&&     \\%
&-c_{(q-1)l+2}&c_{(q-1)l+1}&&&  \\%
&&\ddots&\ddots&& \\%
&&&-c_{ql-1}&c_{ql-1}&  \\%
&&&&-c_{ql}&c_{ql}%
\end{array}%
\right]_{l\times(l+1)},%
\]
for $q=1, 2, \ldots, p-1$,
\[
R_q = \left[ %
\begin{array}{llllll}%
a_{(q-1)l+1}&-a_{(q-1)l+1}&&&&     \\%
&a_{(q-1)l+2}&-a_{(q-1)l+2}&&&  \\%
&&\ddots&\ddots&& \\%
&&&a_{ql-1}&-a_{ql-1}&  \\%
&&&&a_{ql}&-a_{ql}%
\end{array}%
\right]_{l\times(l+1)},%%
\]
\[
L_1 = \left[ %
\begin{array}{llll}%
c_1&&&     \\%
-c_{2}&\ddots&&  \\%
&\ddots&c_{l-1}&  \\%
&&-c_{l}&c_{l}%
\end{array}%
\right]_{l\times l},%
R_p = \left[ %
\begin{array}{llll}%
a_{n-l+1}&-a_{n-l+1}&&     \\%
&a_{n-l+2}&\ddots&  \\%
&&\ddots&-a_{n-1}  \\%
&&&a_{n}%
\end{array}%
\right]_{l\times l},%%
\]
Note that matrices $L_q$ and $R_q$ should be assembled in $G_L$ and
$G_R$ such that their positive diagonal match to the main diagonal
of $G_L$ and $G_R$. Assume that eigenvalues of matrices $G_L$ and
$G_R$ are denoted by vectors $\Lambda_L = [\lambda^l_1, \lambda^l_2,
\ldots, \lambda^l_n]^T$ and $\Lambda_R = [\lambda^r_1, \lambda^r_2,
\ldots, \lambda^r_n]^T$ respectively. The following Lemma determines
$\Lambda_L$ and $\Lambda_R$ explicitly. \\

\begin{lemma} %
\label{lemm:1d:1}%
$G_L$ and $G_R$ are strictly positive definite and we have,
\[
\Lambda_L = [ %
c_1, \ldots, c_l, a_{l+1}, \ldots, a_{2l}, \ldots, %
c_{(p-2)l+1}, \ldots, c_{(p-1)l}, a_{n-l+1}, \ldots, a_{n}%
]^T %
\]%
\[%
\Lambda_R = [ %
a_1, \ldots, a_l, c_{l+1}, \ldots, c_{2l}, \ldots, %
a_{(p-2)l+1}, \ldots, a_{(p-1)l}, c_{n-l+1}, \ldots, c_{n}%
]^T %
\]
\end{lemma}%

\begin{proof}%
Considering the structures of $G_L$ and $G_R$, it is easy to verify
that $G_L$ and $G_R$ are block lower and upper triangular matrices.
Note that the definition of both of the block lower triangular and
upper triangular matrix hold for each of $G_L$ and $G_R$. Therefore,
their corresponding eigenvalues are equal to union of their diagonal
blocks's eigenvalues. The Diagonal blocks have a simple bi-diagonal
structure and their eigenvalues are equal to their main diagonals.
Considering the fact that that $a_i, c_i > 0, \ i=1, \ldots, n$
completes the proof.
\end{proof}%
\\

\begin{theorem}%
\label{th:1d:1}%
The presented parallel one-dimensional algorithm is convergent if
$|1-\omega_L| \  |1-\omega_R| < 1.$\\
\end{theorem}%

\begin{proof}%
Considering the decomposition \eqref{eq:ca1d:decomp} of matrix $A$,
with simple algebra, we can write the presented one-dimensional
parallel iterative algorithm in the following matrix form $(k=0,2, \cdots)$,%
\begin{equation}%
\label{eq:ca1d:matiter}%
\begin{array}{l}
\Big( (\mathbf{b}-\omega_L \Lambda_L) I + \omega_L G_L \Big)\
\mathbf{u}^{(k+1)} =
\Big( ((1-\omega_L)\mathbf{b} + \omega_L \Lambda_R) I - \omega_L G_R
\Big)\
\mathbf{u}^{(k)} +\omega_L \mathbf{f} \\\\ %
\Big( (\mathbf{b}-\omega_R \Lambda_R) I + \omega_R G_R \Big)\ \mathbf{u}^{(k+2)} = %
\Big(((1-\omega_R)\mathbf{b} + \omega_R \Lambda_L) I - \omega_R G_L
\Big)\ \mathbf{u}^{(k+1)}+\omega_R \mathbf{f},
\end{array}
\end{equation}%
where $\mathbf{b} \eqdef [b_1, b_2, \ldots, b_n]^T$. From
\eqref{eq:ca1d:matiter} we have,%
\begin{equation}%
\label{eq:ca1d:matiter2}%
\mathbf{u}^{(k+2)}=T \mathbf{u}^{(k)}+ \bar{\mathbf{f}}, \qquad k = 0, 2, \cdots. %
\end{equation}%
in which,
\begin{eqnarray}%
\label{eq:ca1d:t}%
T = %
\Big( (\mathbf{b}-\omega_R \Lambda_R) I + \omega_R G_R \Big)^{-1}
\Big(((1-\omega_R)\mathbf{b} + \omega_R \Lambda_L) I - \omega_R G_L \Big) %
\nonumber\\
\times %
\Big( (\mathbf{b}-\omega_L \Lambda_L) I + \omega_L G_L \Big)^{-1}%
\Big( ((1-\omega_L)\mathbf{b} + \omega_L \Lambda_R) I - \omega_L G_R \Big)%
\end{eqnarray}%
To proof the convergence of iterations \eqref{eq:ca1d:matiter2} it
is remained to show that the spectral radius of iteration matrix $T$
which is denoted by $\rho(T)$ here is strictly less than unity (cf.
\cite{saad2003ims}). Now, let's to define the matrix $\tilde{T}$
which is equivalent to $T$ (has the same spectral radius), %
\begin{equation}%
\label{eq:ca1d:tt}%
\tilde{T}= %
\Big( (\mathbf{b}-\omega_R \Lambda_R) I + \omega_R G_R \Big)%
 \ T \ %
\Big( (\mathbf{b}-\omega_R \Lambda_R) I + \omega_R G_R \Big)^{-1}%
\end{equation}%
Using \eqref{eq:ca1d:t}, $\tilde{T}$ can be written as follows,%
\begin{eqnarray}%
\label{eq:ca1d:t2}%
\tilde{T}= %
\Big(((1-\omega_R)\mathbf{b} + \omega_R \Lambda_L) I - \omega_R G_L \Big) %
\Big( (\mathbf{b}-\omega_L \Lambda_L) I + \omega_L G_L \Big)^{-1}%
\nonumber\\
\times %
\Big( ((1-\omega_L)\mathbf{b} + \omega_L \Lambda_R) I - \omega_L G_R \Big)%
\Big( (\mathbf{b}-\omega_R \Lambda_R) I + \omega_R G_R \Big)^{-1}
\end{eqnarray}%
Let's to define the following matrices %
\[G_{L-} = ((1-\omega_R)\mathbf{b} + \omega_R \Lambda_L) I - \omega_R G_L, \qquad %
G_{L+} =  (\mathbf{b}-\omega_L \Lambda_L) I + \omega_L G_L \] %
\[G_{R-} = ((1-\omega_L)\mathbf{b} + \omega_L \Lambda_R) I - \omega_L G_R, \qquad %
G_{R+} =  (\mathbf{b}-\omega_R \Lambda_R) I + \omega_R G_R\] %
Assume eigenvalues of $G_{L-}$, $G_{L+}$, $G_{R-}$ and $G_{R+}$ are
denoted by the following vectors respectively, %
\[\Lambda_{L-} = [\lambda_1^{l-}, \ldots, \lambda_n^{l-}]^T, \quad%
\Lambda_{L+} = [\lambda_1^{l+}, \ldots, \lambda_n^{l+}]^T, \] %
\[\Lambda_{R-} = [\lambda_1^{r-}, \ldots, \lambda_n^{r-}]^T, \quad  %
\Lambda_{R+} = [\lambda_1^{r+}, \ldots, \lambda_n^{r+}]^T,\] %
It is easy to show that,%
\[
\Lambda_{L-} = (1-\omega_R) \mathbf{b}, \quad %
\Lambda_{L+} = \mathbf{b}, \quad %
\Lambda_{R-} = (1-\omega_L) \mathbf{b}, \quad %
\Lambda_{R+} = \mathbf{b}.%
\]
Using the knowledge of linear algebra we have, %
\[
\rho(T) = \rho(\tilde{T}) = %
\| \ \tilde{T} \ \| \leqslant %
\| \ (G_{L-}) \ (G_{L+}^{-1}) \ \| \times  %
\| \ (G_{R-}) \ (G_{R+}^{-1}) \ \| %
\]
where the norm operator $\|\cdot\|$ is understood as the Euclidean
norm here,
\[%
\| \ (G_{L-}) \ (G_{L+}^{-1}) \ \| = \max_{1\leqslant i \leqslant n} \ %
\bigg| \ \frac{\lambda_i^{l-}}{\lambda_i^{l+}}\bigg| = %
  \ \max_{1\leqslant i \leqslant n} \ %
\bigg| \ \frac{(1-\omega_R) b_i}{b_i}\bigg| = |1-\omega_R| %
\]%
\[%
\| \ (G_{R-}) \ (G_{R+}^{-1}) \ \| = \max_{1\leqslant i \leqslant n} \ %
\bigg| \ \frac{\lambda_i^{r-}}{\lambda_i^{r+}}\bigg| = %
  \ \max_{1\leqslant i \leqslant n} \ %
\bigg| \ \frac{(1-\omega_L) b_i}{b_i}\bigg| = |1-\omega_L|%
\]
Therefore, if,%
\[
|1-\omega_L| \ |1-\omega_R| < 1.
\]
$\rho(T) < 1$ which complete the proof.
\end{proof}%
\\

\begin{corollary}
The presented one-dimensional parallel Gauss-Seidel algorithm
($\omega_L=\omega_R =
1$) is unconditionally convergent.\\
\end{corollary}

\begin{corollary}
When $\omega_L=\omega_R= \omega$ the convergence criterion of the
presented one-dimensional parallel SOR algorithm is similar to that
of sequential SOR method, i.e., $0<\omega<2$.
Moreover $\omega_{opt} = 1$ in this case.\\
\end{corollary}

\begin{remark}
Using the same procedure, it will easy to show that the coefficient
matrix $A$ corresponding to a high order methods (includes more
extended computational molecules) admits a decomposition like to
\eqref{eq:ca1d:decomp}, where $G_L$ and $G_R$ are two strictly
positive definite block lower and upper triangular matrices
respectively. Then, it will be
straightforward to proof the convergence similarly.\\
\end{remark}

%===============================================================
%===============================================================

\subsection{Two-dimensional convergence analysis}%
\label{sec:conv_anal_2d}

Following the previous section, it is possible to proof the
convergence of the presented method in two spatial dimensions. To
illustrate this issue, the convergence of the presented
two-dimensional algorithm for discrtized equation \eqref{eq:d2d}
will be discussed in this section.

Extension of mono-dimensional proof to higher dimensions is based on
 tensor product properties of the structured grids. The standard five-point
 discretization of \eqref{eq:2d} problem on an $n \times
 m$ structured grid (ignoring Dirichlet boundary nodes)
  leads to the following system of linear
 equation (the effects of Dirichlet boundary nodes are
 included in the right hand side vector),
 \[
 A \mathbf{u} = \mathbf{f}
 \]
 where strictly positive definite $nm \times nm$ dimensional matrix
 $A$ has the following structure (row-wise ordering is assumed for grid points),%
 \[
 A = \left[ %
\begin{array}{ccccc}%
B_1 & A_1 &&& \\ %
C_2 & B_2 & A_2 &&  \\ %
& C_3 & \ddots & \ddots &  \\ %
&& \ddots & \ddots & A_{m-1}   \\ %
&&& C_m  & B_m %
\end{array}%
\right]_{m\times m},
 \]
 where each of $B_i$, $A_i$ and $C_i$ is an $n\times n $
 multi-diagonal (three or mono-diagonal) matrices. For $j=2, \ldots, m-1$,
\[
 B_j = \left[ %
\begin{array}{cccccc}%
a^W_{1,j} & a^P_{1,j} & a^E_{1,j} &&& \\ %
& a^W_{2,j} & \ddots & \ddots &&  \\ %
& & \ddots  & \ddots & a^E_{n-1,j} &  \\ %
&&& a^W_{n,j}  & a^P_{n,j}& a^E_{n,j} %
\end{array}%
\right]_{n\times (n+2)}
\]

\[
 B_1 = \left[ %
\begin{array}{cccccc}%
a^P_{1,1} & a^E_{1,1} &&&& \\ %
a^W_{2,1} & a^P_{2,1} & a^E_{2,1} &&&  \\ %
& a^W_{3,1} & \ddots & \ddots &&  \\ %
&& \ddots & \ddots & a^E_{n-1,1} &   \\ %
&&& a^W_{n,1}  & a^P_{n,1} & a^E_{n,1}%
\end{array}%
\right]_{n\times (n+1)}
\]

\[
 B_m = \left[ %
\begin{array}{cccccc}%
a^W_{1,m} &a^P_{1,m} & a^E_{1,m} &&& \\ %
& a^W_{2,m} & \ddots  & \ddots  &&  \\ %
&& \ddots & \ddots & a^E_{n-2,m} &   \\ %
&&&a^W_{n-1,m} &a^P_{n-1,m} & a^E_{n-1,m} \\ %
&&&& a^W_{n,m}  & a^P_{n,m} %
\end{array}%
\right]_{n\times (n+1)}
\]
and $ A_j = \big[\mathtt{diag}(-a^N_{1,j}, \ldots,
-a^N_{n,j})\big]_{n\times n}, \ C_j = \big[\mathtt{diag}(-a^S_{1,j},
\ldots,
-a^S_{n,j})\big]_{n\times n}.$\\

Without loss of generality, assume that the spatial domain is
decomposed into $p_x \times p_y$ number of sub-domains, where $p_x$
and $p_y$ are even and every sub-domain includes $l_x \times l_y$
number of grid points ($n=p_xl_x$ and $m=p_yl_y$). Also, assume that
the sweeping direction of the first sub-domain is $NE$-to-$SW$ at
the starting iteration. It is easy to show that the coefficient
matrix $A$ admits the following decomposition (notice that the
natural row-wise ordering is considered here),
\begin{equation}%
\label{eq:ca2d:decomp}
 A =  G_1 + \beta I + G_2,
\end{equation}
where $I$ denotes an $n \times m$ identity matrix here,

\[
G_1 = \left[ %
\begin{array}{ccccc}%
N_1&&&&  \\%
&S_2&&&  \\%
&&\ddots&& \\%
&&&N_{p_y-1}&  \\%
&&&&S_{p_y}  %
\end{array}%
\right],\quad%
G_2 = \left[ %
\begin{array}{ccccc}%
S_1&&&&  \\%
&N_2&&&  \\%
&&\ddots&& \\%
&&&S_{p_y-1}&  \\%
&&&&N_{p_y}  %
\end{array}%
\right],%
\]
for $q = 1, \ldots, p_y-1$,
\[
N_q = \left[ %
\begin{array}{cccccc}%
H_{(q-1)l_y+1}&P_{(q-1)l_y+1}&&&&     \\%
&H_{(q-1)l_y+2}&P_{(q-1)l_y+2}&&&     \\%
&&\ddots&\ddots&& \\%
&&&\ddots&P_{ql_y-l}& \\%
&&&&H_{ql_y}&P_{ql_y} %
\end{array}%
\right]_{l_y\times (l_y+1)},%
\]
for $q = 2, \ldots, p_y$,
\[
S_q = \left[ %
\begin{array}{cccccc}%
Q_{(q-1)l_y+1}&V_{(q-1)l_y+1}&&&&     \\%
&Q_{(q-1)l_y+2}&V_{(q-1)l_y+2}&&&     \\%
&&\ddots&\ddots&& \\%
&&&Q_{ql_y-1}&\ddots& \\%
&&&&Q_{ql_y}&V_{ql_y} %
\end{array}%
\right]_{l_y\times (l_y+1)},%
\]
\[
N_{p_y} = \left[ %
\begin{array}{cccc}%
H_{(p_y-1)l_y+1}&P_{(p_y-1)l_y+1}&&     \\%
&\ddots&\ddots& \\%
&&H_{p_yl_y-1}&P_{p_yl_y-l} \\%
&&&H_{p_yl_y} %
\end{array}%
\right]_{l_y\times l_y},%
\]
\[
S_{1} = \left[ %
\begin{array}{cccc}%
V_{1}&&& \\%
Q_{2}&V_{2}& &\\%
&\ddots &\ddots & \\
&&Q_{l_y}&V_{l_y} %
\end{array}%
\right]_{l_y\times l_y},%
\]
where $P_j = A_j$ (for $j=1, \ldots, m-1$), $Q_j=C_j$ (for $j=2,
\ldots, m$) and for $j=1, \ldots, m$,%
\[
H_j = \left[ %
\begin{array}{ccccc}%
NE_{1,j}&&&&  \\%
&NW_{2,j}&&&  \\%
&&\ddots&& \\%
&&&NE_{p_x-1,j}&  \\%
&&&&NW_{p_x,j}  %
\end{array}%
\right],%
\]
\[
V_j = \left[ %
\begin{array}{ccccc}%
SW_{1,j}&&&&  \\%
&SE_{2,j}&&&  \\%
&&\ddots&& \\%
&&&SW_{p_x-1, j}&  \\%
&&&&SE_{p_x,j}  %
\end{array}%
\right],\quad%
\]
denoting $a^{NE}_{i,j} = a^E_{i,j} + a^N_{i,j}$, for $r=1, \ldots,
p_x-1$, $j=1, \ldots m$,
\[
 NE_{r,j} = \left[ %
\begin{array}{cccccc}%
a^{NE}_{(r-1)l_x+1,j} & -a^E_{(r-1)l_x+1,j} &&&& \\\\ %
& a^{NE}_{(r-1)l_x+2,j} & -a^{E}_{(r-1)l_x+2,j} &&& \\\\ %
&& \ddots & \ddots && \\\\ %
&&& a^{NE}_{rl_x-1,j} & -a^{E}_{rl_x-1,j} & \\\\ %
&&&& a^{NE}_{rl_x,j} & -a^E_{rl_x,j}  %
\end{array}%
\right]_{l_x\times (l_x+1)}
\]
denoting $a^{NW}_{i,j} = a^W_{i,j} + a^N_{i,j}$, for $r=2, \ldots,
p_x$, $j=1, \ldots m$,
\[
 NW_{r,j} = \left[ %
\begin{array}{cccccc}%
-a^{W}_{(r-1)l_x+1,j} & a^{NW}_{(r-1)l_x+1,j} &&&& \\\\ %
& -a^{W}_{(r-1)l_x+2,j} & a^{NW}_{(r-1)l_x+2,j} &&& \\\\ %
&& \ddots & \ddots && \\\\ %
&&& -a^{W}_{rl_x-1,j} & a^{NW}_{rl_x-1,j} & \\\\ %
&&&& -a^{W}_{rl_x,j} & a^{NW}_{rl_x,j}  %
\end{array}%
\right]_{l_x\times (l_x+1)}
\]
denoting $a^{SW}_{i,j} = a^W_{i,j} + a^S_{i,j}$, for $r=2, \ldots,
p_x-1$, $j=1, \ldots m$,
\[
 SW_{r,j} = \left[ %
\begin{array}{cccccc}%
-a^{W}_{(r-1)l_x+1,j} & a^{SW}_{(r-1)l_x+1,j} &&&& \\\\ %
& -a^{W}_{(r-1)l_x+2,j} & a^{SW}_{(r-1)l_x+2,j} &&& \\\\ %
&& \ddots & \ddots && \\\\ %
&&& -a^{W}_{rl_x-1,j} & a^{SW}_{rl_x-1,j} & \\\\ %
&&&& -a^{W}_{rl_x,j} & a^{SW}_{rl_x,j}  %
\end{array}%
\right]_{l_x\times (l_x+1)}
\]
denoting $a^{SE}_{i,j} = a^E_{i,j} + a^S_{i,j}$, for $r=2, \ldots,
p_x-1$, $j=1, \ldots m$,
\[
 SE_{r,j} = \left[ %
\begin{array}{cccccc}%
a^{SE}_{(r-1)l_x+1,j} & -a^E_{(r-1)l_x+1,j} &&&& \\\\ %
& a^{SE}_{(r-1)l_x+2,j} & -a^{E}_{(r-1)l_x+2,j} &&& \\\\ %
&& \ddots & \ddots && \\\\ %
&&& a^{SE}_{rl_x-1,j} & -a^{E}_{rl_x-1,j} & \\\\ %
&&&& a^{SE}_{rl_x,j} & -a^E_{rl_x,j}  %
\end{array}%
\right]_{l_x\times (l_x+1)}
\]
Considering effects of boundaries,
\[
 SW_{1,j} = \left[ %
\begin{array}{ccccc}%
 a^{SW}_{1,j} &&& \\\\ %
-a^{W}_{2,j} & a^{SW}_{2,j} && \\\\ %
& \ddots & \ddots & \\\\ %
&& -a^{W}_{l_x-1,j} & a^{SW}_{l_x-1,j} & \\\\ %
&&& -a^{W}_{l_x,j} & a^{SW}_{l_x,j}  %
\end{array}%
\right]_{l_x\times l_x}
\]
\[
 SE_{p_x,j} = \left[ %
\begin{array}{ccccc}%
a^{SE}_{n-l_x+1,j} & -a^E_{(n-l_x+1,j} &&& \\\\ %
& a^{SE}_{n-l_x+2,j} & -a^{E}_{n-l_x+2,j} && \\\\ %
&& \ddots & \ddots & \\\\ %
&&& a^{SE}_{n-1,j} & -a^{E}_{n-1,j}  \\\\ %
&&&& a^{SE}_{n,j}   %
\end{array}%
\right]_{l_x\times l_x}
\]
Note that matrices $N_q$ and $S_q$ should be assembled in $G_1$ and
$G_2$ such that $H_j$ and $V_j$ entries are matched on the
corresponding main diagonals. in the same manner, matrices
$NE_{r,j}$, $NW_{r,j}$, $SE_{r,j}$ and $SW_{r,j}$ should be
assembled in $H_j$ and $V_j$ such that their positive diagonals are
matched on the corresponding main diagonals.

Using the same approach, parallel sweeping along the other diagonal,
$NW$-to-$SE$, implies the following decomposition,
\begin{equation}%
\label{eq:ca2d:decomp2}
 A =  G_3 + \beta I + G_4,
\end{equation}
Since there is a clear analogy between decompositions
\eqref{eq:ca2d:decomp} and \eqref{eq:ca2d:decomp2}, we avoid further
details about $G_3$ and $G_4$  too save the space. In fact, during a
$NE$-to-$SW$ sweep, matrices $G_1$ and $G_2$ respectively include
the implicit and explicit parts of the computational molecule
related to the discrtized diffusion operator. During a $SW$-to-$NE$
sweep, the position of implicit and explicit parts will be reversed.

To simplify the proof of convergence, we define some helpful
matrices and mention their elegant properties here.\\

\begin{definition} Multilevel matrix splitting: Consider $n\times n$ matrix $A$ and
integer set $\Xi = \{m_1, m_2, \ldots, m_s\}$ such that $1 \leqslant
s \leqslant n$, $m_i \geqslant 1$ ($i=1, \ldots, s$) and
$\sum_{i=1}^s m_i = n$. The multilevel splitting of matrix $A$ under
set $\Xi$ which is denoted by $MLS_\Xi(A)$ is defined as set $\{A_1,
A_2, \ldots, A_s\}$ such that every matrix $A_i$ is composed from
the first $m_i$ row and columns of matrix $A_{i-1}$ $(i=1, \ldots,
s)$, where $A_0 = A$.\\
\end{definition}

\begin{definition} Alternatively lower-upper triangular matrix:
An $n\times n$ matrix $A$ is called alternatively lower-upper
triangular, if there is a multilevel splitting set $\Xi = \{m_1,
\ldots, m_s\}$ such that every $A_i \in MLS_\Xi(A)$ is either a
lower triangular matrix or an upper triangular matrix.\\
\end{definition}

\begin{definition} Alternatively block lower-upper triangular matrix:
An $n\times n$ matrix $A$ is called alternatively block lower-upper
triangular, if there is a multilevel splitting set $\Xi = \{m_1,
\ldots, m_s\}$ such that every $A_i \in MLS_\Xi(A)$ is either a
block lower triangular matrix or a block upper triangular matrix.\\
\end{definition}

\begin{corollary}%
\label{coro:2d:1}%
Assume that $n\times n$ matrix $A$ is an alternatively lower-upper
triangular matrix, then $det(A) = \prod_{i=1}^n a_{ii}$ and the set
eigenvalues of matrix $A$ are equal to set of diagonal entries.
\end{corollary}%
\begin{proof}
The proof is evidently followed considering the definition of the
matrix determinant.
\end{proof}\\

\begin{corollary}%
\label{coro:2d:2}%
Assume that $n\times n$ matrix $A$ is an alternatively block
lower-upper triangular matrix, then $det(A)$ is equal to the product
of determinant of its diagonal blocks and the set eigenvalues of
matrix $A$ are equal to union of set of eigenvalues corresponding to
the diagonal blocks.\\
\end{corollary}%

\begin{lemma}%
\label{lemma:2d:3}%
Assume that $n\times n$ matrix $A$ is an alternatively lower-upper
triangular matrix, then there is a set of permutations with
cardinality $m \leqslant n$ that maps $A$ to either a lower
triangular or an upper triangular matrix.\\
\end{lemma}%

\begin{lemma}%
\label{lemma:2d:4}%
Assume that $n\times n$ matrix $A$ is an alternatively block
lower-upper triangular matrix, then there is a set of permutations
with cardinality $m \leqslant n$ that maps $A$ to either a lower
block triangular an upper triangular triangular matrix.
\end{lemma}%
\\
Since Lemma \ref{lemma:2d:3} and \ref{lemma:2d:4} do not play any
role in the course of our analysis, we leave their proof to
interested readers.\\

\begin{corollary}%
\label{coro:2d:5}%
Matrices $G_1$, $G_2$, $G_3$ and $G_4$ in \eqref{eq:ca2d:decomp} and
\eqref{eq:ca2d:decomp2} are alternatively block lower-upper
triangular.
\end{corollary}%
\begin{proof}%
The proof is evident, considering the structures of $G_1$, $G_2$,
$G_3$ and $G_4$ illustrated in this section.
\end{proof}%
\\
It is worth mentioning that the same result is hold for sub-matrices
$H_j$ and $V_j$ which is summarized in the following corollary.\\
\begin{corollary}%
\label{coro:2d:555}%
Matrices $H_j$ and $V_j$ ($j=1, \ldots, m$), are alternatively block
lower-upper triangular.\\
\end{corollary}%

\begin{figure}[ht]%
\centering{%
\includegraphics[width=13.cm]{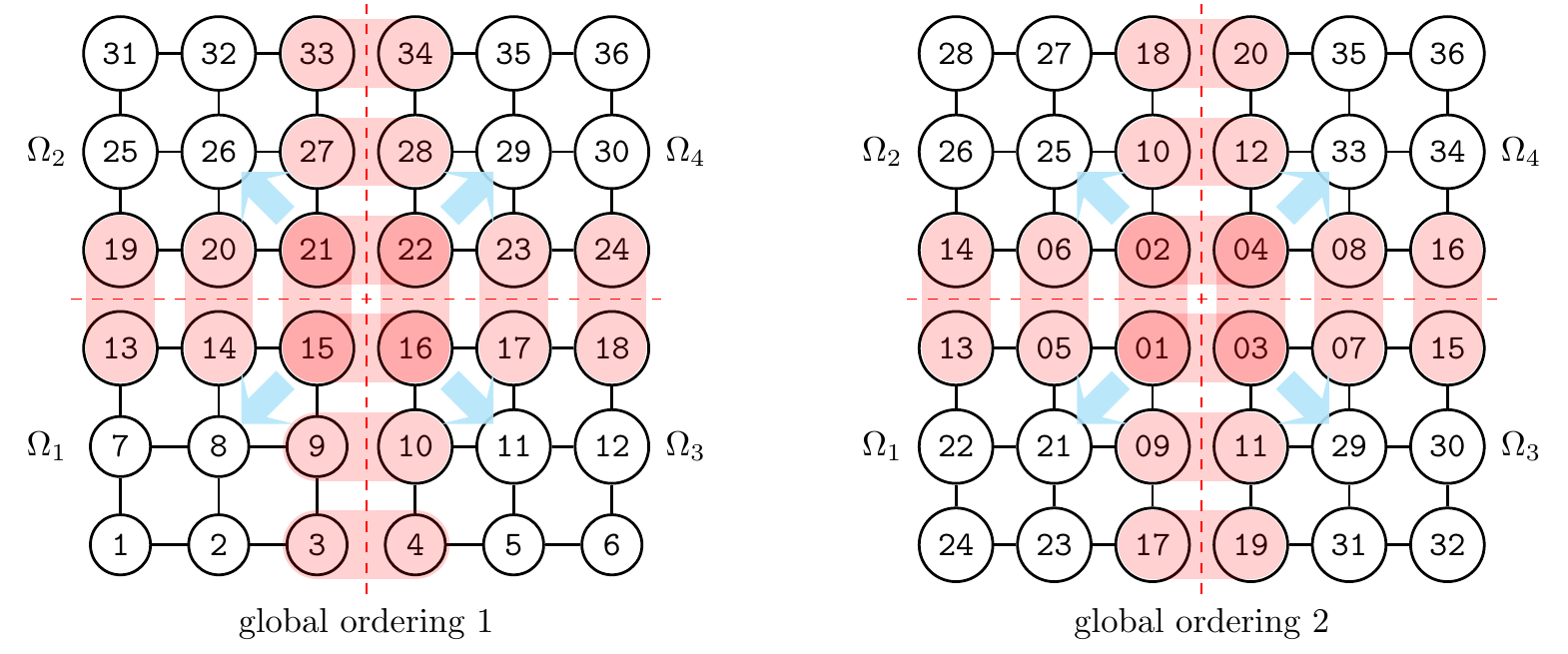}%
}%
\caption{Decomposition of a $6\times6$ Cartesian mesh into
$2\times2$ array of sub-domains: sweeping directions and global
order of updating for natural row-wise ordering (left) and parallel
 sweeping ordering (right).}%
\label{fig:global_ordering}%
\end{figure}%

\begin{figure}[ht]%
\centering{%
\includegraphics[width=13.cm]{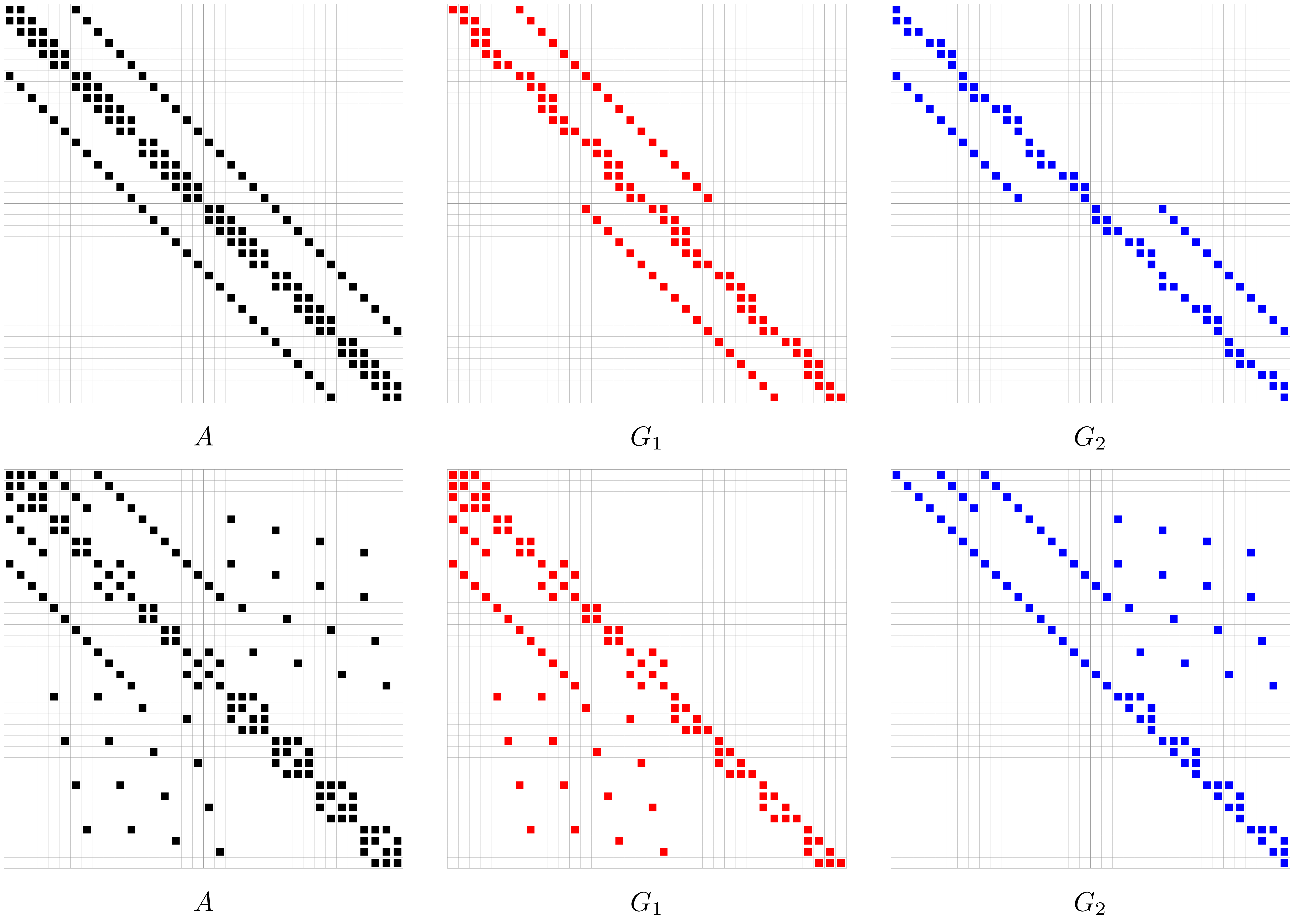}%
}%
\caption{Decomposition of a $6\times6$ Cartesian mesh into
$2\times2$ array of sub-domains: global sparsity pattern of $A$,
$G_1$ and $G_2$ for global ordering based on
 natural row-wise ordering (up) and parallel sweeping ordering (bottom); cf.
  figure \ref{fig:global_ordering}.}%
\label{fig:sparsity_2d}%
\end{figure}

Lets to visualize results of Lemma \ref{lemma:2d:4} and Corollary
\ref{coro:2d:5} through a simple example. For this purpose we use
$2\times2$ decomposition of the $6\times6$ spatial grid shown in
left side of figure \ref{fig:split2d}. To show the global sparsity
pattern, we should first convert local indexing into a global one.
We use two type of global ordering here. The first one based on
natural row-wise ordering (cf. left hand side of figure
\ref{fig:global_ordering}). The second one is globalization of
update order according to our algorithm. Assuming local updating
index of a node \#n is in $\Omega_j$ is denoted by $i_l$, its global
index is computed as $i_g = 4*(i_l-1)+ j$ for coupling nodes and the
remaining nodes are indexed based on priority of their domain ID
(cf. right hand side of figure \ref{fig:global_ordering}). The
global sparsity pattern of $A$, $G_1$ and $G_2$ for these two case
of global ordering is shown in figure \ref{fig:sparsity_2d}. The
plot clearly confirms the results of corollary \ref{coro:2d:5}. It
is worth mentioning that, the global ordering based on the presented
parallel sweeping algorithm defines a permutation to map the
alternatively block lower-upper triangular matrices $G_1$ and $G_2$
to their corresponding block lower (upper) triangular (compare upper
and lower part of figure \ref{fig:sparsity_2d}).

Assume that eigenvalues of matrices $G_1$, $G_2$, $G_3$ and $G_4$ in
\eqref{eq:ca2d:decomp} and \eqref{eq:ca2d:decomp2} are denoted by
 the following vectors vectors respectively %
\[\Lambda_1 = [\lambda^1_1, \lambda^1_2, \ldots, \lambda^1_n]^T, \quad %
\Lambda_2 = [\lambda^2_1, \lambda^2_2, \ldots, \lambda^2_n]^T, \]  %
\[\Lambda_3 = [\lambda^3_1, \lambda^3_2, \ldots, \lambda^3_n]^T, \quad%
\Lambda_4 = [\lambda^4_1, \lambda^4_2, \ldots, \lambda^4_n]^T.\] %
The following Lemma computes these values explicitly.\\
\begin{lemma}%
\label{lemma:2d:6}%
The following relation holds for sets of eigenvalues of $G_1$,
$G_2$, $G_3$ and $G_4$ in \eqref{eq:ca2d:decomp} and
\eqref{eq:ca2d:decomp2},
\[
\Lambda_1,\ \Lambda_2,\ %
\Lambda_3,\ \Lambda_4\ %
\subset \{ %
\ a^{NE}_{i,j},\ %
 a^{NW}_{i,j},\ %
 a^{SE}_{i,j},\ %
 a^{SW}_{i,j}%
\ | \ i=1, \ldots, n; \ j=1, \ldots, m \ \}.
\]
\end{lemma}%
\begin{proof}%
Considering Corollary \ref{coro:2d:5}, the eigenvalues of $G_1$,
$G_2$, $G_3$ and $G_4$ are equal to their diagonal entries. The
diagonal entries are equal to diagonals of $H_i$ and $V_j$ for $j=1,
\ldots, m$ which complete the proof.
\end{proof}%
\\

\begin{corollary}%
\label{coro:2d:7}%
Matrices $G_1$, $G_2$, $G_3$ and $G_4$ in \eqref{eq:ca2d:decomp} and
\eqref{eq:ca2d:decomp2} are strictly positive definite.\\
\end{corollary}%

\begin{theorem}%
\label{th:2d:1}%
The presented parallel two-dimensional algorithm is convergent if
the following condition holds,
\[|1-\omega_{SW}| \ | 1-\omega_{SE}| \ |1-\omega_{NW}| \  |1-\omega_{NE}| < 1.\]
\end{theorem}%

\begin{proof}%
Let's to define the following temporary notations to concise
expressions,
\[\omega_1 \eqdef \omega_{NE}, \quad %
  \omega_2 \eqdef \omega_{SW}, \quad %
  \omega_3 \eqdef \omega_{NW}, \quad %
  \omega_4 \eqdef \omega_{SE}.\]
Inspiring from the 1D analysis, we consider every four iterations of
the presented method as a group of iterations (sweeping along
$NE$-to-$SW$, $SW$-to-$NE$, $NW$-to-$SE$ and $SE$-to-$NW$
respectively). Considering the decomposition \eqref{eq:ca2d:decomp}
and \eqref{eq:ca2d:decomp2} of $A$, we can write the presented
two-dimensional parallel iterative algorithm in the following
matrix form $(k=0,4, \cdots)$,%
\begin{equation}%
\label{eq:ca2d:matiter}%
\begin{array}{l}
\Big( (\mathbf{a}-\omega_1 \Lambda_1) I + \omega_1 G_1 \Big)\
\mathbf{u}^{(k+1)} = \Big( ((1-\omega_1)\mathbf{a} + \omega_1
\Lambda_2) I - \omega_1 G_2 \Big)\
\mathbf{u}^{(k)} \quad + \omega_1\mathbf{f}, \\\\ %
\Big( (\mathbf{a}-\omega_2 \Lambda_2) I + \omega_2 G_2 \Big)\ \mathbf{u}^{(k+2)} = %
\Big(((1-\omega_2)\mathbf{a} + \omega_2 \Lambda_1) I - \omega_2 G_1
\Big)\ \mathbf{u}^{(k+1)}+\omega_2\mathbf{f}, \\\\ %
\Big( (\mathbf{a}-\omega_3 \Lambda_3) I + \omega_3 G_3 \Big)\
\mathbf{u}^{(k+3)} = \Big( ((1-\omega_3)\mathbf{a} + \omega_3
\Lambda_4) I - \omega_3 G_4 \Big)\
\mathbf{u}^{(k+2)} + \omega_3\mathbf{f}, \\\\ %
\Big( (\mathbf{a}-\omega_4 \Lambda_4) I + \omega_4 G_4 \Big)\ \mathbf{u}^{(k+4)} = %
\Big(((1-\omega_4)\mathbf{a} + \omega_4 \Lambda_3) I - \omega_4 G_3
\Big)\ \mathbf{u}^{(k+3)}+\omega_4\mathbf{f}, \\\\ %
\end{array}
\end{equation}%
where vector $\mathbf{a}$ includes diagonal entries of $A$ with an
order consistent to $G_i$. From
\eqref{eq:ca2d:matiter} we have,%
\begin{equation}%
\label{eq:ca2d:matiter2}%
\mathbf{u}^{(k+4)}=T \mathbf{u}^{(k)}+ \bar{\mathbf{f}}, \qquad k=0,4, \cdots. %
\end{equation}%
where,
\begin{eqnarray}%
\label{eq:ca2d:t}%
T = %
\Big((\mathbf{a}-\omega_4 \Lambda_4) I + \omega_4 G_4 \Big)^{-1}%
\Big(((1-\omega_4)\mathbf{a}+\omega_4\Lambda_3)I-\omega_4 G_3 \Big)%
\nonumber\\
\times %
\Big( (\mathbf{a}-\omega_3 \Lambda_3) I + \omega_3 G_3 \Big)^{-1}%
\Big(((1-\omega_3)\mathbf{a}+\omega_3\Lambda_4)I-\omega_3 G_4 \Big)%
\nonumber\\
\times %
\Big((\mathbf{a}-\omega_2 \Lambda_2) I + \omega_2 G_2 \Big)^{-1}%
\Big(((1-\omega_2)\mathbf{a}+\omega_2 \Lambda_1)I-\omega_2 G_1 \Big)%
\nonumber\\
\times %
\Big( (\mathbf{a}-\omega_1 \Lambda_1) I + \omega_1 G_1 \Big)^{-1}%
\Big(((1-\omega_1)\mathbf{a}+\omega_1\Lambda_2)I-\omega_1 G_2 \Big)%
\end{eqnarray}%
With some simple manipulations of iteration matrix $T$, we have,%
\begin{equation}%
\label{eq:ca2d:tt}%
\tilde{T} =  %
(G_{1-})\ (G_{1+})^{-1}\ %
(G_{2-})\ (G_{2+})^{-1}\ %
(G_{3-})\ (G_{3+})^{-1}\ %
(G_{4-})\ (G_{4+})^{-1} %
\end{equation}%
such that $\rho(T) = \rho(\tilde{T})$, and,
%
%1
\[%
G_{1+} = (\mathbf{a}-\omega_1 \Lambda_1) I + \omega_1 G_1, \qquad %
G_{1-} = ((1-\omega_2)\mathbf{a}+\omega_2 \Lambda_1)I-\omega_2 G_1 %
\]%
%
% 2
\[%
G_{2+} = (\mathbf{a}-\omega_2 \Lambda_2) I + \omega_2 G_2, \qquad %
G_{2-} = ((1-\omega_1)\mathbf{a}+\omega_1\Lambda_2)I-\omega_1 G_2 %
\]%
%
% 3
\[%
G_{3+} = (\mathbf{a}-\omega_3 \Lambda_3) I + \omega_3 G_3, \qquad %
G_{3-} = ((1-\omega_4)\mathbf{a}+\omega_4\Lambda_3)I-\omega_4 G_3 %
\]%
%
% 4
\[%
G_{4+} = (\mathbf{a}-\omega_4 \Lambda_4) I + \omega_4 G_4, \qquad %
G_{4-} = ((1-\omega_3)\mathbf{a}+\omega_3\Lambda_4)I-\omega_3 G_4 %
\]%
Assume eigenvalues of $G_{1-}$, $G_{1+}$, $G_{2-}$, $G_{2+}$,
$G_{3-}$, $G_{3+}$, $G_{4-}$ and $G_{4+}$ are denoted by vectors
$\Lambda_{1-} = \{\lambda_i^{1-}\}$, $\Lambda_{1+} =
\{\lambda_i^{1+}\}$, $\Lambda_{2-} = \{\lambda_i^{1-}\}$,
$\Lambda_{2+} = \{\lambda_i^{2+}\}$, $\Lambda_{3-} =
\{\lambda_i^{3-}\}$, $\Lambda_{3+} = \{\lambda_i^{3+}\}$,
$\Lambda_{4-} = \{\lambda_i^{4-}\}$ and $\Lambda_{4+} =
\{\lambda_i^{4+}\}$ (for $i=1, \ldots, mn$) respectively,
it is easy to show that,%
\[
\Lambda_{1-} = (1-\omega_2) \mathbf{a}, \quad %
\Lambda_{1+} = \mathbf{a}, \quad %
\Lambda_{2-} = (1-\omega_1) \mathbf{a}, \quad %
\Lambda_{2+} = \mathbf{a},%
\]
\[
\Lambda_{3-} = (1-\omega_4) \mathbf{a}, \quad %
\Lambda_{3+} = \mathbf{a}, \quad %
\Lambda_{4-} = (1-\omega_3) \mathbf{a}, \quad %
\Lambda_{4+} = \mathbf{a}.%
\]
Therefore, we have ($\|\cdot\|$ denotes the Euclidean norm),%
\begin{align}%
\label{eq:2d:tts}%
\nonumber
\| \tilde{T} \| & \leqslant %
\| (G_{1-})\ (G_{1+})^{-1} \| \ %
\| (G_{2-})\ (G_{2+})^{-1} \| \ %
\| (G_{3-})\ (G_{3+})^{-1} \| \  %
\| (G_{4-})\ (G_{4+})^{-1}\| \\ \nonumber %
& =
\Bigg(%
\max_{1\leqslant i \leqslant n} \ %
\bigg| \ \frac{\lambda_i^{1-}}{\lambda_i^{1+}}\bigg| %
\Bigg) \
\Bigg(%
\max_{1\leqslant i \leqslant n} \ %
\bigg| \ \frac{\lambda_i^{2-}}{\lambda_i^{2+}}\bigg|%
\Bigg) \ %
\Bigg(%
\max_{1\leqslant i \leqslant n} \ %
\bigg| \ \frac{\lambda_i^{4-}}{\lambda_i^{3+}}\bigg| %
\Bigg) \ %
\Bigg(%
\max_{1\leqslant i \leqslant n} \ %
\bigg| \ \frac{\lambda_i^{4-}}{\lambda_i^{4+}}\bigg|%
\Bigg) %
\\ & =
|1-\omega_2|\ |1-\omega_1|\ |1-\omega_4|\ |1-\omega_3|.
\end{align}
Since, $\rho(T) =  \rho(\tilde{T}) = \| \tilde{T} \|$, using
\eqref{eq:2d:tts}, the following condition ensures the convergence
of iterations: $|1-\omega_1|\ |1-\omega_2|\ |1-\omega_3|\
 |1-\omega_4| < 1.$
\end{proof}%
\\

\begin{corollary}
The presented two-dimensional parallel Gauss-Seidel algorithm
($\omega_{NE}=\omega_{SW}=\omega_{NW}=\omega_{SE}=1$) is unconditionally convergent.\\
\end{corollary}

\begin{corollary}
When $\omega_{NE}=\omega_{SW}=\omega_{NW}=\omega_{SE}=\omega$ the
convergence criterion of the presented two-dimensional parallel SOR
algorithm is similar to that of sequential SOR method, i.e.,
$0<\omega<2$. Moreover, in this case $\omega_{opt} = 1$.\\
\end{corollary}

%===============================================================
%===============================================================

\section{Extension to general structured n-diagonal matrices}%
\label{sec:gextension}

Using the Cartesian tensor product properties of structured grids
and following the same procedure mentioned in the previous sections,
the extension of our method to three (or higher) spatial dimensions
is straightforward. Therefore, it is not of value to further discuss
on this issue. Now, let us to generalize our method a bit more.\\

\begin{definition} (Cartesian grid) %
The spatial network $Z \subset \mathbb{R}^d$ $(d\geqslant 1)$ is
called as a Cartesian grid in $\mathbb{R}^d$ if it is formed by
tensor product of $d$ one-dimensional spatial grids. \\ %
\end{definition}

\begin{definition} (Structured grid) %
The spatial network $Z \subset \mathbb{R}^d$ $(d\geqslant 1)$ is
called as a structured grid in $\mathbb{R}^d$ if there is an
orientation preserving deformation in $\mathbb{R}^d$ which maps $Z$
to a Cartesian grid in $\mathbb{R}^d$. \\ %
\end{definition}

It is clear that every structured grid in $\mathbb{R}^d$ has $2^d$
corner points and so $2^d$ directed diagonal directions (sweeping
directions in our algorithm).\\

\begin{definition} (Structured n-diagonal matrix)
\label{def:snm} %
Consider the following assumptions: $Z \subset \mathbb{R}^d$ is
structured grid with $m_i$ grid points ($i=1, \ldots, d$) along each
spatial dimension. Every grid point of $Z$ has an ID number based on
a natural row-wise ordering. $A$ is a square matrix with cardinality
$n = \prod_{i=1}^d m_i$. $A$ is an n-diagonal matrix with a
symmetric sparsity pattern ($n$ is odd). Then we called $A$ as a
structured n-diagonal matrix if there is a computational molecule
$M$ which maps nonzero entries of every row $i$ of $A$ to a set of
points in the neighborhood of  node $i$ in $Z$; and vice versa.
 \\ %
\end{definition}
We assume that the reader is aware from the required special
treatment at/near boundaries of $Z$. Therefore, these issues are
ignore within Definition \ref{def:snm}. In fact by Definition
\ref{def:snm} we would like to extend application of our method to
every matrix $A$ which can be connected to a structured grid, e.g.,
it can be raised from discretization of a PDE on a structured grid.

Now, we want to generalize previous results, i.e., solving
(preconditioning) system of linear equations $A \mathbf{u} =
\mathbf{f}$ by parallel SOR (or ILU) when $A$ is an irreducible
diagonal dominant structured n-diagonal matrix. Let us to briefly
mention the outline of algorithm in this case. Assume $Z \subset
\mathbb{R}^d$ is the corresponding structured grid to $A$. We first
decompose $Z$ into parts (based on load balancing criteria) and
define the sweeping directions for each sub-domain during the course
of each $2^d$ iterations. Then, order of update for coupling points
are determined; note that the coupling points are determined based
on the current sweeping direction and the computational molecule.
Finally the iterative procedure is started by updating coupling
points in the corresponding order and then internal nodes. Now let
us to re-state the convergence theories for this general case.
Although, we do not proof the following results, their proof should
be straightforward (though can be very massive if one does not find
an appropriate abstraction in this regard).\\

\begin{theorem}%
\label{th:6:gen}%
Consider irreducible diagonal dominant structured n-diagonal matrix
$A$ with the corresponding structured grid $Z \subset \mathbb{R}^d$.
The the following results hold.%
\begin{itemize}
\item [(i)] Matrix admit $(2^d)/2$ decomposition in the following form
such that every decomposition correspond to parallel sweeping along
an appropriate diagonal of $Z$: $A = G_{j} + D + G_{j+1}$ for $j=1,
3, \ldots, d-1$, where $D$ is a diagonal matrix with nonnegative
entries, for $i=1, \ldots, d$, each $G_{i}$ is strictly positive
definite alternatively block lower-upper triangular matrix with
eigenvalues equal to main diagonal.
\item [(ii)] Assume that for $i=1, \ldots, d$; $\omega_i$ denotes
the relaxation parameter along the $i$-th sweeping direction. Then
the presented parallel SOR algorithm is convergent under the
following condition: $\prod_{i=1}^{2^d} |1-\omega_i| <1$.
\end{itemize}
\end{theorem}

%===============================================================
%===============================================================

\section{Numerical experiment}%
\label{sec:res}%

\noindent In this section we present numerical results on the
convergence and performance of the presented method. As the
computing platform we used a 32-core mechine based on 16 Dual-Core
Intel 2400 Mhz nodes (4GB RAM on each node) switched with a Gbit
Ethernet network. The MPICH \cite{bridges1995ugm} library was used
for message passing implementation of exchange algorithm.

The $d$-dimensional ($d=1, 2, 3$) version of following laplace
equation is used here as a model problem,
\begin{equation}%
\label{eq:numer:model}%
\sum_{i=1}^{i=d} \bigg(\frac{\partial^2 u}{\partial x_i^2}\bigg) = 0 %
 \quad \mathtt{in} \quad [0,1]^d, %
\qquad u = \prod_{i=1}^{i=d} x_i \quad \mathtt{on} \quad \partial
\big([0,1]^d \big)
\end{equation}%
where $x_i$ denotes the $i$-th component of spatial position vector,
e.g., for $d=3$, $\mathbf{x} = (x_1, x_2, x_3) = (x, y, z)$. The
exact solution of \eqref{eq:numer:model} is $\hat{u}(\mathbf{x}) =
\prod_{i=1}^{i=d} x_i$. Since this exact solution is spatially
$d$-linear, it is possible to recover the exact solution by a
spatially second order numerical method, independent of the grid
resolution. Therefore, the discretization errors does can not
contribute in our numerical analysis. In all of our numerical
examples in this study, the $L_1$-norm of error (with respect to the
exact solution) is used to study the convergence. Assume that the
total number of degree of freedoms is $n$ after the discretization
of \eqref{eq:numer:model}, the $L_1$-norm of error is computed as
follows, $L_1-\mathtt{error} = \frac{1}{n}\ {\sum_{l=1}^{n} | u_l -
\hat{u}_l|}$, %
where $u_i$ and $\hat{u}_i$ denote respectively the approximate and
exact solutions at the $i$-th degree of freedom. Moreover, the
initial guess is taken to be $u_l = 0$ (for $l=1, \ldots, n$) and
the iterations is stopped when the $L_1$-norm of error is below
threshold $1.e-3$; except in our 3D experiments in which the
convergence threshold $1.e-2$ is taken into account.

Although we do not plan to release a software related to the
presented method, some (crude) Fortran codes related to model
problems used in this study are freely available on the following
web (mainly to demonstrate the validity of results presented in this
section):
$\href{http://sites.google.com/site/rohtav/home/par}%
{http://sites.google.com/site/rohtav/home/par}.$

\subsection{One-dimensional numerical experiment}

In our one-dimensional numerical experiments \eqref{eq:numer:model}
is solved on uniform grids with resolutions, 41, 81 and 161.

In the first 1D example, the model problem is solved by In this left
to right sweeping Gauss-Seidel (LRGS), right to left sweeping
Gauss-Seidel (RLGS), symmetric sweeping Gauss-Seidel (SGS) and the
presented parallel Gauss-Seidel with various number of sub-domains
(PGS(p), p denotes the number of sub-domains).

Table \ref{tab:1d:1} shows the results of this experiment. It is
clear that the presented method is convergent independent from the
number of sub-domains, and the number of iterations to meet the
convergence criterion is close to that of the original Gauss-Seidel
methods.

%%%%%%%%%%%%%%%%%%%%%%%%%%%%%%%%%%%%%%%%%%%%%%%%%%%%%%%%%%%%%%%
%%%%%%%%%%%%%%%%%%%%%%%%%%%%%%%%%%%%%%%%%%%%%%%%%%%%%%%%%%%%%%%

\begin{table}[ht]%
\caption{1D parallel Gauss-Seidel vs. 1D sequential Gauss-Seidel.}%
\label{tab:1d:1}%
\footnotesize\rm
\begin{tabular*}{\textwidth}{@{}r l*{7}{@{\extracolsep{0pt plus12pt}}l}}
\hline%
&&41 grids&&81 grids&&161 grids&\\
\cline{3-4} \cline{5-6} \cline{7-8} \\
method&
&iteration&$L_1$-error&iteration&$L_1$-error&iteration&$L_1$-error\\\hline%
LRGS   &&979&9.94266e-4&3905&9.98916e-4&15598&9.99738e-4\\
RLGS   &&960&9.94266e-4&3866&9.98916e-4&15519&9.99738e-4\\
SGS    &&976&9.96647e-4&3892&9.99752e-4&15565&9.99977e-4\\
PGS(2) &&975&9.95198e-4&3891&9.99297e-4&15563&9.99825e-4\\
PGS(4) &&974&9.97789e-4&3890&9.99986e-4&15561&9.99779e-4\\
PGS(6) &&974&9.94523e-4&3890&9.99231e-4&15559&9.99769e-4\\
PGS(8) &&973&9.97436e-4&3890&9.98472e-4&15557&9.99807e-4\\
PGS(10) &&973&9.94144e-4&3889&9.99221e-4&15555&9.99809e-4\\
PGS(14) &&972&9.93927e-4&3888&9.99208e-4&15551&9.99755e-4\\
PGS(18) &&970&9.99382e-4&3887&9.99179e-4&15547&9.99783e-4\\
PGS(24) &&-&-&3885&9.99857e-4&15541&9.99814e-4\\
PGS(30) &&-&-&3884&9.99459e-4&15535&9.99738e-4\\
PGS(36) &&-&-&3882&9.99864e-4&15529&9.99704e-4\\
\hline
\end{tabular*}
\end{table}

%%%%%%%%%%%%%%%%%%%%%%%%%%%%%%%%%%%%%%%%%%%%%%%%%%%%%%%%%%%%%%%
%%%%%%%%%%%%%%%%%%%%%%%%%%%%%%%%%%%%%%%%%%%%%%%%%%%%%%%%%%%%%%%

In the second 1D example the model problem is solved by left to
right sweeping SOR (LRSOR), right to left sweeping SOR (RLSOR),
symmetric sweeping SOR (SSOR) and the presented parallel SOR with
various number of sub-domains (PSOR(p), p denotes the number of
sub-domains). The optimal values of $\omega$ are applied in these
experiments. These values are primarily computed via numerical
experiments (searching in $[1.0,2.0]$ with $1.e-3$ accuracy).

Tables \ref{tab:1d:2}, \ref{tab:1d:3} and \ref{tab:1d:4} show
results of these experiments. According to tables, PSOR method is
convergent and has a competitive convergence rate with sequential
SOR methods. Since the initial error had an asymmetric spatial
distribution, the sweeping directions has considerable effect on the
iteration counts. For example RLSOR method has the best convergence
rate. Note the presented parallel solver overcomes LRSOR and SSOR in
some cases.

%%%%%%%%%%%%%%%%%%%%%%%%%%%%%%%%%%%%%%%%%%%%%%%%%%%%%%%%%%%%%%%
%%%%%%%%%%%%%%%%%%%%%%%%%%%%%%%%%%%%%%%%%%%%%%%%%%%%%%%%%%%%%%%

\begin{table}[ht]
\caption{1D parallel SOR vs. 1D sequential SOR for grid resolution
41.}%
\label{tab:1d:2}%
\footnotesize\rm %
\begin{tabular*}{\textwidth}{@{}r l*{6}{@{\extracolsep{0pt plus12pt}}l}}
\hline%
method& &iteration&$L_1$-error&$\omega_{L}$- optimum&$\omega_{R}$-
optimum\\\hline%
LRSOR &&51& 9.99298e-4& 1.86887&    -\\
RLSOR &&31&    9.38294e-4& -      &    1.86637\\
SSOR &&62    &9.48613e-4 &1.00000    &1.87776\\
PSOR(2) &&31    &9.46469e-4 &1.84970    &1.92084\\
PSOR(4) &&76    &9.88365e-4 &1.00000    &1.94890\\
PSOR(6) &&78    &9.88246e-4 &1.00000    &1.94890\\
PSOR(8) &&72    &9.85974e-4 &1.00000    &1.89379\\
PSOR(10) &&71   &9.72100e-4 &1.00000    &1.88577\\
PSOR(14) &&87   &9.59733e-4 &1.00000    &1.95591\\
PSOR(18) &&71   &9.60855e-4 &1.00000    &1.90080\\
\hline
\end{tabular*}
\end{table}

%%%%%%%%%%%%%%%%%%%%%%%%%%%%%%%%%%%%%%%%%%%%%%%%%%%%%%%%%%%%%%%
%%%%%%%%%%%%%%%%%%%%%%%%%%%%%%%%%%%%%%%%%%%%%%%%%%%%%%%%%%%%%%%

\begin{table}[ht]
\caption{1D parallel SOR vs. 1D sequential SOR for grid resolution 81.}%
\label{tab:1d:3}%
\footnotesize\rm%
\begin{tabular*}{\textwidth}{@{}r l*{6}{@{\extracolsep{0pt plus12pt}}l}}
\hline method& &iteration&$L_1$-error&$\omega_{L}$-
optimum&$\omega_{R}$- optimum\\ \hline%
LRSOR &&103&    9.82824e-4& 1.93193&    -      \\
RLSOR &&61&    9.42462e-4& -      &    1.93143\\
SSOR &&120&  8.87046e-4& 1.00000&    1.93487 \\
PSOR(2) &&80&   9.54235e-4& 1.90982&    1.95691\\
PSOR(4) &&164&  9.93129e-4& 1.00000&    1.97194 \\
PSOR(6) &&129&  9.66998e-4& 1.00000&    1.96092 \\
PSOR(8) &&141&  9.91414e-4& 1.00000&    1.94589 \\
PSOR(10) &&140& 9.80610e-4& 1.00000&    1.94088 \\
PSOR(14) &&129& 9.35418e-4& 1.00000&    1.95992\\
PSOR(18) &&138& 9.96443e-4& 1.00000&    1.95090 \\
PSOR(24) &&141& 9.93632e-4& 1.00000&    1.94790 \\
PSOR(30) &&180& 9.91055e-4& 1.00000&    1.97094 \\
PSOR(36) &&143& 9.91698e-4& 1.00000&    1.94389 \\%
 \hline
\end{tabular*}
\end{table}

%%%%%%%%%%%%%%%%%%%%%%%%%%%%%%%%%%%%%%%%%%%%%%%%%%%%%%%%%%%%%%%
%%%%%%%%%%%%%%%%%%%%%%%%%%%%%%%%%%%%%%%%%%%%%%%%%%%%%%%%%%%%%%%

%%%%%%%%%%%%%%%%%%%%%%%%%%%%%%%%%%%%%%%%%%%%%%%%%%%%%%%%%%%%%%%
%%%%%%%%%%%%%%%%%%%%%%%%%%%%%%%%%%%%%%%%%%%%%%%%%%%%%%%%%%%%%%%

\begin{table}[ht]
\caption{1D parallel SOR vs. 1D sequential SOR for grid resolution
161.}%
\label{tab:1d:4}%
\footnotesize\rm %
\begin{tabular*}{\textwidth}{@{}r l*{6}{@{\extracolsep{0pt plus12pt}}l}}
\hline%
method& &iteration&$L_1$-error&$\omega_{L}$- optimum&$\omega_{R}$-
optimum\\\hline%
LRSOR &&208&    9.81709e-4& 1.96593&    -      \\
RLSOR &&122&   9.58420e-4& -      &    1.96493\\
SSOR &&236&  9.99534e-4& 1.19840&    1.96693\\
PSOR(2) &&233&  9.82315e-4& 1.00000&    1.96593\\
PSOR(4) &&342&  9.84557e-4& 1.00000&    1.98497\\
PSOR(6) &&253&  9.85584e-4& 1.00000&    1.97996\\
PSOR(8) &&279&  9.72206e-4& 1.00000&    1.97194\\
PSOR(10) &&277& 9.88102e-4& 1.00000&    1.96994\\
PSOR(14) &&280& 9.98979e-4& 1.00000&    1.96994\\
PSOR(18) &&269& 9.91725e-4& 1.00000&    1.97495\\
PSOR(24) &&277& 9.89229e-4& 1.00000&    1.97395\\
PSOR(30) &&278& 9.81003e-4& 1.00000&    1.96894\\
PSOR(36) &&281& 9.73789e-4& 1.00000&    1.97194\\%
\hline
\end{tabular*}
\end{table}

%%%%%%%%%%%%%%%%%%%%%%%%%%%%%%%%%%%%%%%%%%%%%%%%%%%%%%%%%%%%%%%
%%%%%%%%%%%%%%%%%%%%%%%%%%%%%%%%%%%%%%%%%%%%%%%%%%%%%%%%%%%%%%%

Since the optimal values of $\omega_{L}$ in SSOR and PSOR solvers
are close to unity, it suggests to apply simultaneous
over-relaxation and under-relaxation to improve the convergence.
Results of symmetric successive over/under-relaxation SOR (SSOUR)
and parallel successive over/under-relaxation SOR (PSOUR) are
compared in Table \ref{tab:1d:5} for grid resolution $41$ (the
optimal values of relaxation factors are searched in $[0.0,2.0]$
with accuracy $1.e-3$).

%%%%%%%%%%%%%%%%%%%%%%%%%%%%%%%%%%%%%%%%%%%%%%%%%%%%%%%%%%%%%%%
%%%%%%%%%%%%%%%%%%%%%%%%%%%%%%%%%%%%%%%%%%%%%%%%%%%%%%%%%%%%%%%

\begin{table}
\caption{1D parallel over/under-relaxation method vs. 1D
sequential SOR for grid resolution 41.}%
\label{tab:1d:5}%
\footnotesize\rm%
\begin{tabular*}{\textwidth}{@{}r l*{6}{@{\extracolsep{0pt plus12pt}}l}}
\hline%
method& &iteration&$L_1$-error&$\omega_{L}$- optimum&$\omega_{R}$-
optimum \\\hline%
LRSOR &&51& 9.99298e-4& 1.86887&    -      \\
RLSOR &&31&    9.38294e-4& -      &    1.86637\\
SSOUR &&58&  9.99474e-4& 0.24825&    1.87087\\
PSOUR(2) &&33&  8.68232e-4& 1.84785&    1.91892\\
PSOUR(4) &&57&  9.99032e-4& 0.17317&    1.87588\\
PSOUR(6) &&57&  9.98090e-4& 0.11712&    1.87387\\
PSOUR(8) &&58&  9.99746e-4& 0.09910&    1.87087\\
PSOUR(10) &&57& 9.97719e-4& 0.18118&    1.87287\\
PSOUR(14) &&52& 9.99428e-4& 0.36837&    1.89590\\
PSOUR(18) &&53& 9.98391e-4& 0.33233&    1.87988\\%
\hline
\end{tabular*}
\end{table}

%%%%%%%%%%%%%%%%%%%%%%%%%%%%%%%%%%%%%%%%%%%%%%%%%%%%%%%%%%%%%%%
%%%%%%%%%%%%%%%%%%%%%%%%%%%%%%%%%%%%%%%%%%%%%%%%%%%%%%%%%%%%%%%

%===============================================================
%===============================================================

\subsection{Two-dimensional numerical experiment}%

In our two-dimensional numerical experiments \eqref{eq:numer:model}
is solved by row-wise ordering Gauss-Seidel (RGS) and SOR (RSOR),
symmetric Gauss-Seidel (SGS) and SOR (SSOR), frontal sweeping
Gauss-Seidel (FGS) and SOR (FSOR) and the presented parallel
Gauss-Seidel (PGS($p_x$, $p_y$), where $p_x$ and $p_y$ are the
number of sub-domains in $x$ and $y$ directions respectively) and
SOR (PSOR($p_x$, $p_y$)) on grid resolutions $51 \times 51$, $101
\times 101$ and $151 \times 151$.

\begin{figure}[ht]%
\centering{%
\includegraphics[width=13cm]{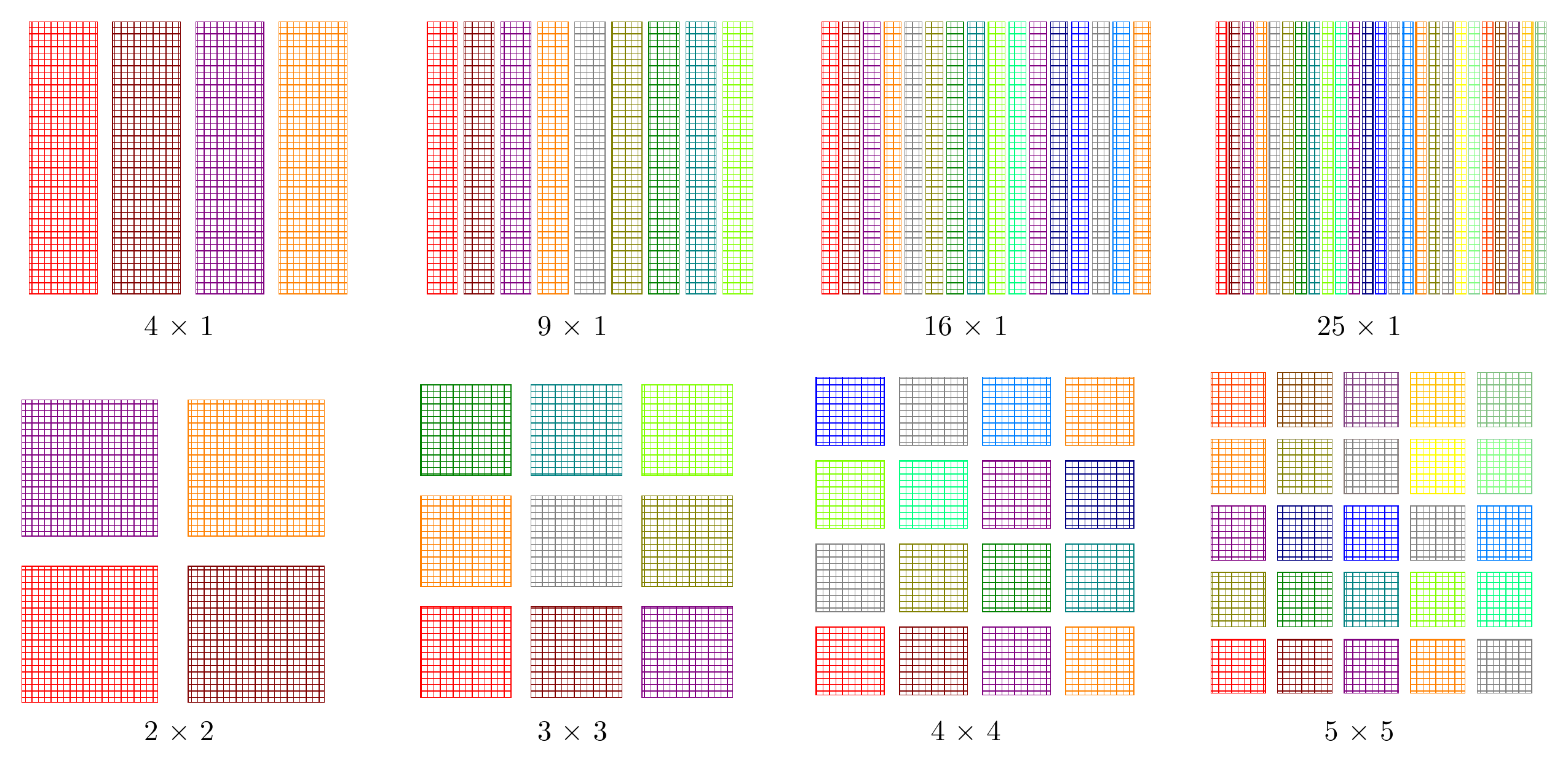}%
}%
\caption{Strip-wise (up) and square-wise (down) sub-domains
topologies used in this study.}%
\label{fig:2ddd}%
\end{figure}

Since the topology of domain decomposition may changes the
convergence rate of parallel solver, strip-wise and square-wise
topologies are examined in this experiment (cf. figure
\ref{fig:2ddd}). The strip-wise and square-wise domain
decompositions have a $p \times 1$ and $\sqrt{p} \times \sqrt{p}$
topologies respectively, where $p$ is the number of sub-domains.

Results of this experiment are summarized in Tables \ref{tab:2d:1},
\ref{tab:2d:2} and \ref{tab:2d:3}. According to tables, the number
of iterations to meet the convergence for the presented parallel
method is close to that of sequential ones, and there is not
sensible dependency to the topology of domain decomposition. Note
that in this experiment the relaxation parameters for all directions
are taken to be equal.

 \clearpage
%%%%%%%%%%%%%%%%%%%%%%%%%%%%%%%%%%%%%%%%%%%%%%%%%%%%%%%%%%%%%%%
% OMEGA =1.0
%%%%%%%%%%%%%%%%%%%%%%%%%%%%%%%%%%%%%%%%%%%%%%%%%%%%%%%%%%%%%%%
\begin{table}
\caption{2D parallel Gauss-Seidel vs. 2D sequential Gauss-Seidel.}%
\label{tab:2d:1}%
\footnotesize\rm
\begin{tabular*}{\textwidth}{@{}l l*{7}{@{\extracolsep{0pt plus12pt}}l}}
\hline%
grid&&51 $\times$ 51 &&101 $\times$ 101&& 151 $\times$ 151&\\
 \cline{3-4} \cline{5-6} \cline{7-8}
method&
&iteration&$L_1$-error&iteration&$L_1$-error&iteration&$L_1$-error\\\hline
 RGS                                                &&         1018&
  9.98560e-4 &         4065 &   9.99123e-4 &         9139 &
  9.99798e-4 \\
 SGS                                                &&         1006&
  9.98720e-4 &         4038 &   9.99374e-4 &         9097 &
  9.99960e-4 \\
 FGS                                                &&         1006&
  9.96308e-4 &         4037 &   9.99753e-4 &         9097 &
  9.99686e-4 \\
 PGS(4$\times$1)                                  &&         1020&
  9.97194e-4 &         4066 &   9.99808e-4 &         9140 &
  9.99839e-4 \\
 PGS(2$\times$2)                                  &&         1020&
  9.97652e-4 &         4065 &   9.99717e-4 &         9138 &
  9.99888e-4 \\
 PGS(9$\times$1)                                  &&         1038&
  9.97011e-4 &         4103 &   9.99232e-4 &         9195 &
  9.99884e-4 \\
 PGS(3$\times$3)                                  &&         1029&
  9.99224e-4 &         4082 &   9.99687e-4 &         9163 &
  9.99873e-4 \\
 PGS(25$\times$1)                                  &&         1088&
  9.96490e-4 &         4219 &   9.99572e-4 &         9371 &
  9.99965e-4 \\
 PGS(5$\times$5)                                  &&         1049&
  9.96392e-4 &         4116 &   9.99972e-4 &         9213 &
  9.99814e-4 \\
\hline
\end{tabular*}
\end{table}

%%%%%%%%%%%%%%%%%%%%%%%%%%%%%%%%%%%%%%%%%%%%%%%%%%%%%%%%%%%%%%%
% OMEGA =1.25
%%%%%%%%%%%%%%%%%%%%%%%%%%%%%%%%%%%%%%%%%%%%%%%%%%%%%%%%%%%%%%%
\begin{table}
\caption{2D parallel SOR vs. 2D sequential SOR. for $\omega=1.25$.}%
\label{tab:2d:2}%
\footnotesize\rm
\begin{tabular*}{\textwidth}{@{}l l*{7}{@{\extracolsep{0pt plus12pt}}l}}
\hline%
grid&&51 $\times$ 51 &&101 $\times$ 101&& 151 $\times$ 151&\\
 \cline{3-4} \cline{5-6} \cline{7-8} %
method&
&iteration&$L_1$-error&iteration&$L_1$-error&iteration&$L_1$-error\\\hline
RSOR                                               && 616&
  9.97982e-4 &         2450 &   9.99198e-4 &         5501 &
  9.99339e-4 \\
SSOR                                               && 606&
  9.95519e-4 &         2425 &   9.98966e-4 &         5461 &
  9.99332e-4 \\
FSOR                                               && 605&
  9.95328e-4 &         2424 &   9.98924e-4 &         5460 &
  9.99310e-4 \\
PSOR(4$\times$1)                                 &&          626&
  9.95846e-4 &         2467 &   9.98854e-4 &         5524 &
  9.99442e-4 \\
PSOR(2$\times$2)                                 &&          626&
  9.98955e-4 &         2465 &   9.99991e-4 &         5521 &
  9.99318e-4 \\
PSOR(9$\times$1)                                 &&          652&
  9.96804e-4 &         2520 &   9.99830e-4 &         5605 &
  9.99658e-4 \\
 PSOR(3$\times$3)                                 &&          637&
  9.99291e-4 &         2487 &   9.99869e-4 &         5556 &
  9.99914e-4 \\
PSOR(16$\times$1)                                 &&          685&
  9.94546e-4 &         2593 &   9.99467e-4 &         5718 &
  9.99331e-4 \\
PSOR(4$\times$4)                                 &&          648&
  9.98578e-4 &         2512 &   9.99153e-4 &         5590 &
  9.99347e-4 \\
PSOR(25$\times$1)                                 &&          715&
  9.94756e-4 &         2688 &   9.99302e-4 &         5863 &
  9.99984e-4 \\
PSOR(5$\times$5)                                 &&          662&
  9.94443e-4 &         2535 &   9.99867e-4 &         5624 &
  9.99977e-4 \\
\hline
\end{tabular*}
\end{table}

%%%%%%%%%%%%%%%%%%%%%%%%%%%%%%%%%%%%%%%%%%%%%%%%%%%%%%%%%%%%%%%
% OMEGA =1.5
%%%%%%%%%%%%%%%%%%%%%%%%%%%%%%%%%%%%%%%%%%%%%%%%%%%%%%%%%%%%%%%
\begin{table}
\caption{2D parallel SOR vs. 2D sequential SOR. for $\omega=1.5$.}%
\label{tab:2d:3}%
\footnotesize\rm %
\begin{tabular*}{\textwidth}{@{}l l*{7}{@{\extracolsep{0pt plus12pt}}l}}%
\hline%
grid&&51 $\times$ 51 &&101 $\times$ 101&& 151 $\times$ 151&\\
 \cline{3-4} \cline{5-6} \cline{7-8} %
method&
&iteration&$L_1$-error&iteration&$L_1$-error&iteration&$L_1$-error\\\hline
RSOR                                               && 348&
  9.90645e-4 &         1373 &   9.98771e-4 &         3074 &
  9.99890e-4 \\
 SSOR                                               &&          341&
  9.88632e-4 &         1351 &   9.99049e-4 &         3038 &
  9.98968e-4 \\
FSOR                                               && 339&
  9.89902e-4 &         1349 &   9.99485e-4 &         3036 &
  9.99170e-4 \\
PSOR(4$\times$1)                                 &&          369&
  9.97722e-4 &         1415 &   9.98874e-4 &         3133 &
  9.99910e-4 \\
PSOR(2$\times$2)                                 &&          369&
  9.89711e-4 &         1410 &   9.98715e-4 &         3127 &
  9.99310e-4 \\
PSOR(9$\times$1)                                 &&          407&
  9.97695e-4 &         1498 &   9.97574e-4 &         3259 &
  9.99389e-4 \\
PSOR(3$\times$3)                                 &&          382&
  9.99520e-4 &         1443 &   9.97418e-4 &         3179 &
  9.98991e-4 \\
PSOR(16$\times$1)                                 &&          453&
  9.94660e-4 &         1606 &   9.96612e-4 &         3432 &
  9.99365e-4 \\
PSOR(4$\times$4)                                 &&          396&
  9.94275e-4 &         1474 &   9.99090e-4 &         3227 &
  9.98820e-4 \\
PSOR(25$\times$1)                                 &&          477&
  9.90157e-4 &         1736 &   9.99730e-4 &         3645 &
  9.99353e-4 \\
PSOR(5$\times$5)                                 &&          407&
  9.99724e-4 &         1504 &   9.98626e-4 &         3274 &
  9.99382e-4 \\%
\hline
\end{tabular*}
\end{table}

%===============================================================
%===============================================================

\clearpage%

\subsection{Thee-dimensional numerical experiment}%

\begin{figure}[ht]%
\centering{%
\includegraphics[width=13cm]{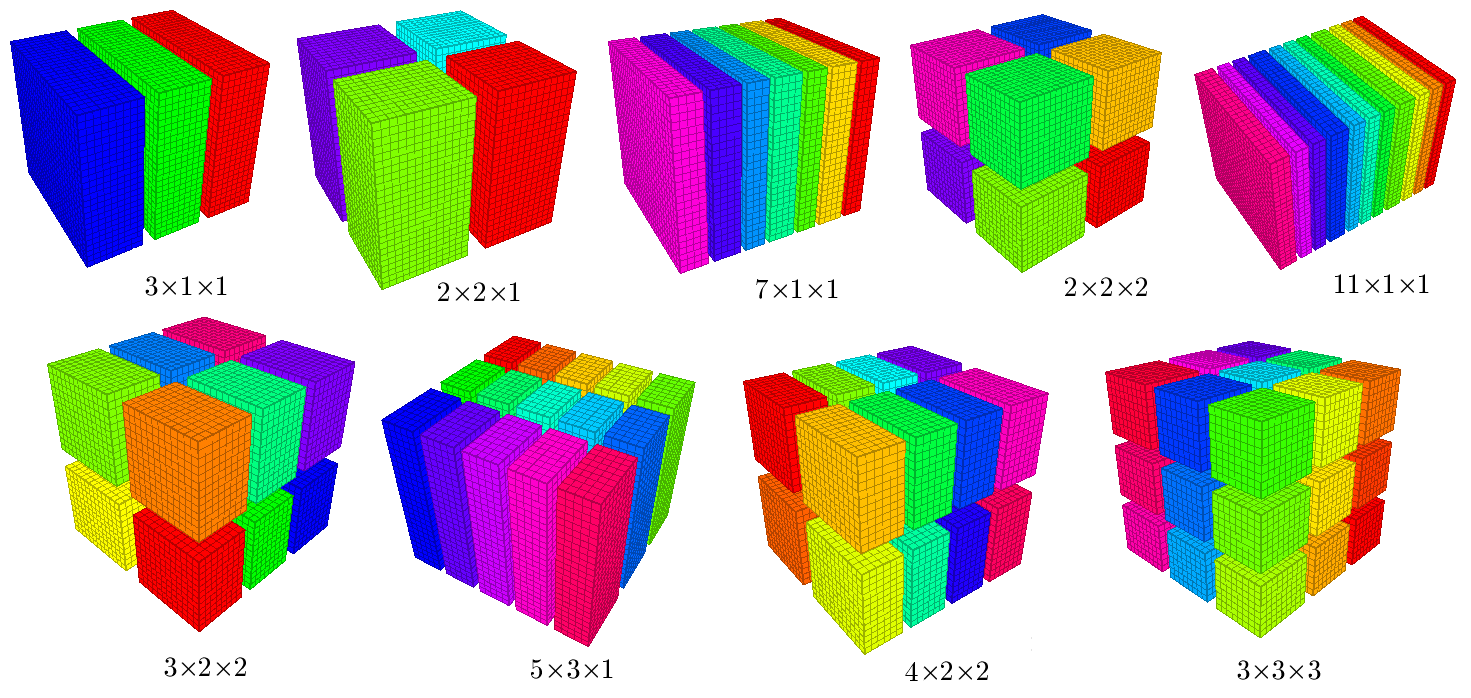}%
}%
\caption{Topologies of domain decompositions in 3D numerical
experiment of the present study.}%
\label{fig:3ddd}%
\end{figure}

In our three-dimensional numerical experiments
\eqref{eq:numer:model} is solved by row-wise sweeping Gauss-Seidel
(RGS) and SOR (RSOR), symmetric sweeping Gauss-Seidel (SGS) and SOR
(SSOR), frontal sweeping Gauss-Seidel (FGS) and SOR (FSOR) and the
presented parallel Gauss-Seidel (PGS($p_x$, $p_y$, $p_z$), where
$p_x$, $p_y$ and $p_z$ are number of sub-domains in x, y and z
directions respectively) and SOR (PSOR($p_x$, $p_y$, $p_z$)) on grid
resolutions $25 \times 25 \times 25$, $51 \times 51 \times 51$ and
$101 \times 101\times 101$. Figure \ref{fig:3ddd} shows the
topologies of domain decompositions used in this numerical
experiment. The relaxation parameters are taken to be equal along
all sweeping directions in this section.

Table \ref{tab:3d:1}, \ref{tab:3d:2} and \ref{tab:3d:3} show results
of this numerical experiment. Tables show that the presented
parallel 3D solvers are convergent and number of iterations for
achieving convergence are close to that of sequential solvers, in
particular in some cases the parallel solver performs slightly
better than the sequential solver. Like 2D results, the convergence
rate of 3D solver has a little dependency to topology of domain
decomposition in this experiment.

In contrast to previously suggested parallel SOR methods in
literature our results are very promising, as we do not observe
sensible convergence decay due to parallelization which is a
pathological weakness of alternative parallel SOR methods. This
observation suggest us to use the presented parallel SOR solver as a
cache efficient SOR method which is disseised in the nest
section.\\\\

%

%%%%%%%%%%%%%%%%%%%%%%%%%%%%%%%%%%%%%%%%%%%%%%%%%%%%%%%%%%%%%%%
% OMEGA =1.0
%%%%%%%%%%%%%%%%%%%%%%%%%%%%%%%%%%%%%%%%%%%%%%%%%%%%%%%%%%%%%%%
\begin{table}[ht]
\caption{3D parallel Gauss-Seidel vs. 3D sequential Gauss-Seidel.}%
\label{tab:3d:1}%
\footnotesize\rm
\begin{tabular*}{\textwidth}{@{}l l*{7}{@{\extracolsep{0pt plus12pt}}l}}
\hline%
grid&&25$\times$25$\times$25 &&51$\times$51$\times$51&&101$\times$101$\times$101&\\
 \cline{3-4} \cline{5-6} \cline{7-8}
method&
&iteration&$L_1$-error&iteration&$L_1$-error&iteration&$L_1$-error\\\hline
RGS                                                && 110&
  9.92077e-3 &          480 &   9.96665e-3 &         1921 &
  9.99432e-3 \\
SGS                                                && 104&
  9.95897e-3 &          466 &   9.99655e-3 &         1893 &
  9.99775e-3 \\
 FGS                                                &&          104&
  9.86877e-3 &          466 &   9.97495e-3 &         1893 &
  9.99243e-3 \\
PGS(3$\times$3$\times$3)                              && 105&
  9.96386e-3 &          468 &   9.99384e-3 &         1898 &
  9.99252e-3 \\
PGS(2$\times$2$\times$1)                              && 105&
  9.99993e-3 &          469 &   9.98302e-3 &         1898 &
  9.99899e-3 \\
PGS(7$\times$1$\times$1)                              && 107&
  9.87250e-3 &          472 &   9.97408e-3 &         1905 &
  9.99024e-3 \\
PGS(2$\times$2$\times$2)                              && 106&
  9.93976e-3 &          470 &   9.99815e-3 &         1901 &
  9.99735e-3 \\
PGS(11$\times$1$\times$1)                              && 108&
  9.92317e-3 &          475 &   9.98559e-3 &         1911 &
  9.99717e-3 \\
PGS(3$\times$2$\times$2)                              && 107&
  9.92017e-3 &          472 &   9.96988e-3 &         1904 &
  9.99496e-3 \\
PGS(5$\times$3$\times$1)                              && 108&
  9.94702e-3 &          474 &   9.97328e-3 &         1907 &
  9.99631e-3 \\
PGS(4$\times$2$\times$2)                              && 108&
  9.83564e-3 &          473 &   9.97272e-3 &         1906 &
  9.99269e-3 \\
PGS(3$\times$3$\times$3)                              && 109&
  9.91718e-3 &          475 &   9.97940e-3 &         1909 &
  9.99410e-3 \\
\hline
\end{tabular*}
\end{table}

%%%%%%%%%%%%%%%%%%%%%%%%%%%%%%%%%%%%%%%%%%%%%%%%%%%%%%%%%%%%%%%
% OMEGA =1.25
%%%%%%%%%%%%%%%%%%%%%%%%%%%%%%%%%%%%%%%%%%%%%%%%%%%%%%%%%%%%%%%
\begin{table}[ht]
\caption{3D parallel SOR vs. 3D sequential SOR. for $\omega=1.25$.}%
\label{tab:3d:2}%
\footnotesize\rm
\begin{tabular*}{\textwidth}{@{}l l*{7}{@{\extracolsep{0pt plus12pt}}l}}
\hline%
grid&&25$\times$25$\times$25 &&51$\times$51$\times$51&&101$\times$101$\times$101&\\
 \cline{3-4} \cline{5-6} \cline{7-8}
method&
&iteration&$L_1$-error&iteration&$L_1$-error&iteration&$L_1$-error\\\hline%
RSOR                                               && 69&
  9.77818e-3 &          293 &   9.99745e-3 &         1164 &
  9.99186e-3 \\
SSOR                                               && 63&
  9.92592e-3 &          281 &   9.93894e-3 &         1137 &
  9.98696e-3 \\
FSOR                                              &&           62&
  9.95549e-3 &          280 &   9.94525e-3 &         1136 &
  9.98851e-3 \\
PGS(3$\times$3$\times$3)                             && 65&
  9.75364e-3 &          284 &   9.97072e-3 &         1144 &
  9.98524e-3 \\
PGS(2$\times$2$\times$1)                            && 65&
  9.87121e-3 &          284 &   9.99866e-3 &         1144 &
  9.99540e-3 \\
PGS(7$\times$1$\times$1)                             && 67&
  9.85868e-3 &          289 &   9.95523e-3 &         1154 &
  9.99055e-3 \\
PGS(2$\times$2$\times$2)                            && 66&
  9.89078e-3 &          287 &   9.93992e-3 &         1149 &
  9.98397e-3 \\
PGS(11$\times$1$\times$1)                             && 69&
  9.73977e-3 &          293 &   9.99839e-3 &         1164 &
  9.99148e-3 \\
PGS(3$\times$2$\times$2)                             && 67&
  9.96697e-3 &          288 &   9.98010e-3 &         1152 &
  9.99076e-3 \\
PGS(5$\times$3$\times$1)                              && 68&
  9.81936e-3 &          291 &   9.95989e-3 &         1157 &
  9.99021e-3 \\
PGS(4$\times$2$\times$2)                             && 68&
  9.84237e-3 &          290 &   9.95890e-3 &         1155 &
  9.99221e-3 \\
PGS(3$\times$3$\times$3)                             && 69&
  9.90408e-3 &          292 &   9.99708e-3 &         1159 &
  9.99217e-3 \\
\hline
\end{tabular*}
\end{table}

%%%%%%%%%%%%%%%%%%%%%%%%%%%%%%%%%%%%%%%%%%%%%%%%%%%%%%%%%%%%%%%
% OMEGA =1.5
%%%%%%%%%%%%%%%%%%%%%%%%%%%%%%%%%%%%%%%%%%%%%%%%%%%%%%%%%%%%%%%
\begin{table}[hb]
\caption{3D parallel SOR vs. 3D sequential SOR. for $\omega=1.5$.}%
\label{tab:3d:3}%
\footnotesize\rm%
\begin{tabular*}{\textwidth}{@{}l l*{7}{@{\extracolsep{0pt plus12pt}}l}}
\hline
grid&&25$\times$25$\times$25 &&51$\times$51$\times$51&&101$\times$101$\times$101&\\
 \cline{3-4} \cline{5-6} \cline{7-8}
method&
&iteration&$L_1$-error&iteration&$L_1$-error&iteration&$L_1$-error\\\hline
RSOR                                               && 41&
  9.82565e-3 &          169 &   9.97880e-3 &          659 &
  9.99566e-3 \\
SSOR                                               && 36&
  9.90569e-3 &          157 &   9.95554e-3 &          633 &
  9.97890e-3 \\
FSOR                                               && 35&
  9.57206e-3 &          155 &   9.99694e-3 &          631 &
  9.98964e-3 \\
PGS(3$\times$3$\times$3)                           && 39&
  9.58539e-3 &          162 &   9.93716e-3 &          643 &
  9.99297e-3 \\
PGS(2$\times$2$\times$1)                              && 38&
  9.81012e-3 &          162 &   9.99555e-3 &          644 &
  9.99660e-3 \\
PGS(7$\times$1$\times$1)                            && 41&
  9.90818e-3 &          170 &   9.90995e-3 &          659 &
  9.99539e-3 \\
PGS(2$\times$2$\times$2)                              && 40&
  9.87981e-3 &          166 &   9.95387e-3 &          651 &
  9.99802e-3 \\
PGS(11$\times$1$\times$1)                              && 44&
  9.61852e-3 &          178 &   9.90053e-3 &          675 &
  9.99266e-3 \\
PGS(3$\times$2$\times$2)                            && 41&
  9.72949e-3 &          169 &   9.93594e-3 &          655 &
  9.98725e-3 \\
PGS(5$\times$3$\times$1)                            && 42&
  9.84760e-3 &          172 &   9.90824e-3 &          662 &
  9.98893e-3 \\
PGS(4$\times$2$\times$2)                             && 42&
  9.66169e-3 &          170 &   9.93732e-3 &          660 &
  9.97495e-3 \\
PGS(3$\times$3$\times$3)                            && 43&
  9.73555e-3 &          173 &   9.92163e-3 &          664 &
  9.99070e-3 \\
\hline
\end{tabular*}
\end{table}

%%%%%%%%%%%%%%%%%%%%%%%%%%%%%%%%%%%%%%%%%%%%%%%%%%%%%%%%%%%%%%%
%%%%%%%%%%%%%%%%%%%%%%%%%%%%%%%%%%%%%%%%%%%%%%%%%%%%%%%%%%%%%%%

\clearpage%

%===============================================================
%===============================================================

\subsection{Cache-efficient SOR relaxation}
\label{sec:chache}

The growing speed gap between processor and memory has led to the
development of hierarchical memories and utility of using caches in
modern processors \cite{sellappa2004cem}. Nowadays the speed of a
code depends increasingly on how well the cache structure is
exploited in the course of computation. The number of cache misses
is recently an important factor parallel to the number of floating
point operations (FLOPS), during comparison of numerical algorithms.
Unfortunately, it is not easy to a-priori estimate number of cache
misses, but it is easy to perform a posteriori analysis
\cite{douglas2000cos}.

When a program references a memory location, the data in the
referenced location and nearby locations are brought into a cache
level. Any additional references to data before these are evicted
from the cache will be one or two orders of magnitude faster than
references to main memory. So increasing the data locality leads to
decreasing cache miss and so the computational performance. A
program is said to have a good data locality if, most of the time,
the computer finds the data referenced by the program in its cache.
The locality of a program can be improved by changing the order of
computation and/or the assignment of data to memory locations so
that references to the same or nearby locations occur near to each
other during the program execution \cite{strout2004sts}. Interested
readers are refereed to
\cite{douglas2000cos,im2004sof,sellappa2004cem,strout2004sts} for
further details about this topic.

Based on the above discussion, each parallel algorithm based on the
domain decomposition concept is potentially a cache efficient
algorithm. The presented method in this study decomposes the global
data into smaller sub-domains (with more locality), and after
getting essential information from neighbor sub-domains at start of
each iteration, performs computations in each sub-domain
independently. In this section we shall study the performance of
presented parallel SOR from the viewpoint of cache efficiency. For
this purpose the 3D model problem is solved with the presented
parallel Gauss-Seidel method ($\omega=1.0$) on grid resolutions $25
\times 25 \times 25$, $51 \times 51 \times 51$, $101 \times 101
\times 101$ and $151 \times 151 \times 151$ (a single processor is
used in this experiment).

\begin{figure}[ht]%
\centering{%
\includegraphics[width=9.cm]{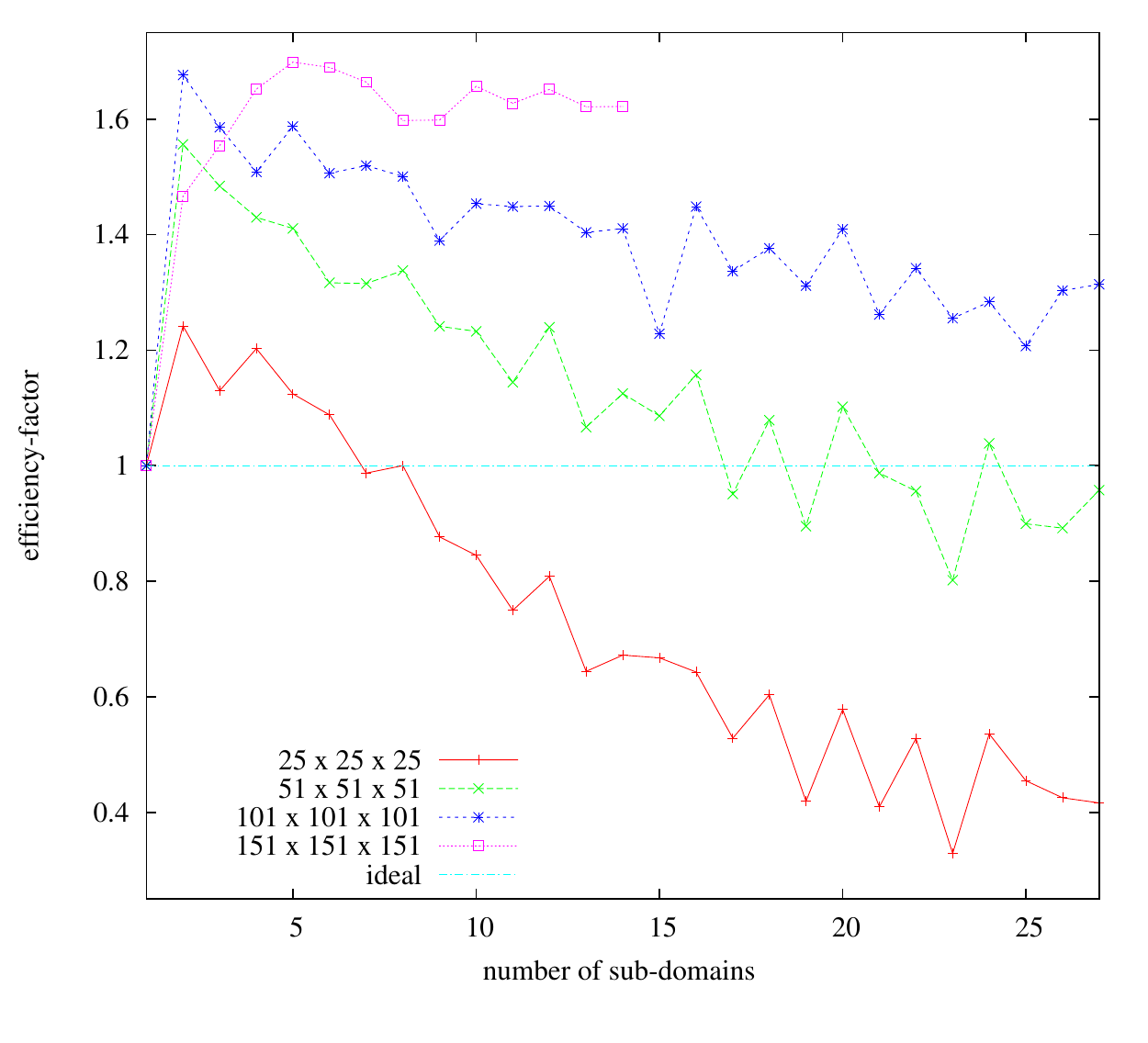}%
}%
\caption{Parallel SOR as a cache-efficient solver: The variation of
efficiency factor vs. number of sub-domains for 3D model problem on
grid resolutions $25^3$, $51^3$, $101^3$ and $151^3$.}%
\label{fig:eff_factor}%
\end{figure}

Let's the define the efficiency-factor as the ratio of the total CPU
time for cache efficient solver to that of classic sequential
solver. Figure \ref{fig:eff_factor} shows variation of
efficiency-factor vs. number of sub-domains in this experiment. The
plot shows that when the size of global data is sufficiently large
(the cache-miss is susceptible), the efficiency-factor is increased
with increasing the number of sub-domains (data locality). Note that
with increasing the number of sub-domain the number of FLOPS is
slightly increased, but due to mentioned effect solver performs
superior.

%

%===============================================================
%===============================================================

\section{Closing remarks}%
\label{sec:closing_remark}%

In the present study, a new method is developed to parallelize
sequential sweeping procedure on structured grids, and is used to
parallelize the SOR (and ILU) preconditioning method on structured
grids. This method can be considered as a special overlapping domain
decomposition which uses a fully parallel multi-frontal sweeping
strategy with local coupling at sub-domains's interfaces. The
implementation of method in one and two spatial dimensions is
discussed in details and its extension to higher dimension is
briefly outlined. The application of method then extended to general
structured n-diagonal matrices. Using notion of alternatively block
upper-lower triangular matrices, the convergence theory is
established in general. Our numerical show that the convergence rate
of the presented method is close to that of sequential solver. It is
indeed a very promising result does not reported for previous
parallel versions of SOR to our knowledge (in fact our method does
not show sensible convergence rate decay due to parallilization).
This result suggest us to use this solver as cache efficient SOR
method. Numerical results on this issue also support effectiveness
of this idea.

It is worth mentioning that our strategy can be directly sued to
parallelize other sequential procedures on structured grid. For
example in the supplementary material to this work (see:
\cite{rtparfa2009}), we use the same strategy to parallelize the
fast sweeping method to solve the eikonal equation. The efficiency
of that method in some cases (complex geometries) are significantly
better than the original fast sweeping as well as classic fast
marching method. Moreover in \cite{rtparsut2005}, we introduce a
fully explicit and unconditionally stable parallel method to solve
nonlinear transient heat equation. This method is based on
parallelizing multi-dimensional version of the Saul'yev scheme using
the presented marching algorithm in this study. Finally, we believe
that the presented strategy can be extended to semi-structured grids
(like quadtree/octree grids). At the moment, we think that under
certain assumptions such an extension is straightforward. Lets to
outline this idea on quadtrees. We denote by Cartesian graded
quadtree a grid which is resulted from an isotropic hierarchical
graded refinement of a square-shaped root cell. Now consider
quadtree grid $\mathcal{Q}$. Assume there is an orientation
preserving deformation $\mathcal{T}$ which maps $\mathcal{Q}$ into a
graded Cartesian quadtree denoted by $\widehat{\mathcal{Q}}$.
Moreover, assume that there is a square-shaped decomposition
\footnote{a decomposition with $\sqrt{p} \times \sqrt{p}$ topology,
where $p$ is number of sub-domains.} $\mathcal{D}$ on
$\widehat{\mathcal{Q}}$ such that every sub-domain is still a graded
Cartesian quadtree. Then it is possible to apply the presented
parallel strategy on $\mathcal{D}$. In this case the sweeping
directions within each sub-domain are determined based on $Z$-order
space filing cures (cf. \cite{berger2005performance}) started from
sub-domain's corners.

%===============================================================
%===============================================================

\section*{Acknowledgements} We would like to thanks triple referees
for their constructive comments which significantly improve the
contents of this paper.

%===============================================================
%===============================================================

\bibliographystyle{plain} % other options: plainnat, plain, abbrv, unsrt, alpha
\bibliography{biblio}%

%===============================================================
%===============================================================

\end{document}